\documentclass[12pt]{amsart}

\usepackage{amscd, graphicx, psfrag}

 \usepackage{ifthen}    
 \usepackage{amssymb}   
 \usepackage{amstext}
 \usepackage{amsfonts}
 \usepackage{amsmath}   
 \usepackage{amscd}     
 \usepackage{latexsym}  
 \usepackage{color}
 \usepackage{epsfig}
 \usepackage{epic}
 \usepackage{pstricks,pst-node,pst-text,pst-tree}
 \usepackage[all]{xy}
 \usepackage{array}
 \usepackage{amsthm}
 \usepackage{eucal, mathrsfs} 
\usepackage{hyperref}

\newcommand{\numberset}[1]{\ensuremath{\mathbb{#1}}}    
\newcommand{\C}{\numberset{C}}  
\newcommand{\N}{\numberset{N}}  
\newcommand{\R}{\numberset{R}}  
\newcommand{\Z}{\numberset{Z}}  
\newcommand{\PP}{\numberset{P}}  

\renewcommand{\L}{\mathcal L}

\newcommand{\p}{\mathfrak p}

\newcommand{\DF}{\mathcal D^bFuk}
\newcommand{\DC}{\mathcal D^bCoh}
\newcommand{\rsh}{R\underline{Hom}}
\renewcommand{\O}{\mathcal O}

\newcommand{\I}{\mathcal I}
\newcommand{\D}{\mathcal D}
\newcommand{\F}{\mathcal F}
\newcommand{\T}{\mathcal T}
\newcommand{\CC}{\mathcal C}
\newcommand{\Cc}{\mathfrak C}

\newcommand{\ris}{\stackrel{\sim}{\longrightarrow}}
\newcommand{\sky}{\mathcal O}

\renewcommand{\t}{\mathfrak t}

\newcommand{\znor}{Z_{\text{nor}}}
\newcommand{\zbnor}{\bar{Z}_{\text{nor}}}
\newcommand{\gnor}{\Gamma_{\text{nor}}}

\theoremstyle{definition}

\newtheorem{thm}{Theorem}[section]
\newtheorem{prop}[thm]{Proposition}
\newtheorem{lem}[thm]{Lemma}
\newtheorem{cor}[thm]{Corollary}
\newtheorem{rem}[thm]{Remark}
\newtheorem{ex}[thm]{Example}
\newtheorem{defi}[thm]{Definition}

\newtheorem{ass}[thm]{Assumption}
\newtheorem{con}[thm]{Conjecture}

\DeclareMathOperator{\im}{Im} 
 \DeclareMathOperator{\Arg}{Arg}
 
\DeclareMathOperator{\Crit}{Crit} \DeclareMathOperator{\supp}{supp}

\DeclareMathOperator{\id}{Id} 
\DeclareMathOperator{\spn}{span} 
\DeclareMathOperator{\inv}{inv} \DeclareMathOperator{\Log}{Log}
\DeclareMathOperator{\rk}{rk}

\newboolean{mynotes}\setboolean{mynotes}{true}

\newcommand{\mycomments}[1]{
           \ifthenelse{\boolean{mynotes}}
                      {#1}{}
           }

\begin{document}

\title{Symmetries of Lagrangian fibrations}
\date{August 7, 2009}
\author{Ricardo Casta\~no-Bernard, Diego Matessi, Jake P. Solomon}

\begin{abstract}
We construct fiber-preserving anti-symplectic involutions for a
large class of symplectic manifolds with Lagrangian torus
fibrations. In particular, we treat the K3 surface and the six
dimensional examples constructed in \cite{CB-M}, which include a six
dimensional symplectic manifold homeomorphic to the quintic
threefold. We interpret our results as corroboration of the view
that in homological mirror symmetry, an anti-symplectic involution
is the mirror of duality. In the same setting, we construct
fiber-preserving symplectomorphisms that can be interpreted as the
mirror to twisting by a holomorphic line bundle.
\end{abstract}

\maketitle

\tableofcontents

\section{Introduction}
\subsection{Statement of result}
Let $X$ be a symplectic manifold and let $B$ be smooth manifold of
half the dimension of $X.$ We call a continuous map $f: X
\rightarrow B$ a \emph{Lagrangian fibration} if each fiber of $f$
contains a relatively open dense set that is a smooth Lagrangian
submanifold of $X.$ Lagrangian fibrations arose classically in the
context of integrable systems and toric geometry. More recently,
Lagrangian fibrations have played a role in the conjectural
interpretation of mirror symmetry introduced by Strominger, Yau and
Zaslow \cite{SYZ}. We discuss this in greater detail in
Section~\ref{ssec:m}.

In \cite{CB-M}, the
first two authors introduced a general construction that produced a
class $\CC$ of Lagrangian fibrations. See Section~\ref{sec:bk} for
the precise definition of $\CC.$ In short, $\CC$ consists of
fibrations on symplectic manifolds of dimensions $4$ and $6.$ In
those dimensions, $\CC$ includes fibrations such that the total
space is homeomorphic to any of the Calabi-Yau complete
intersections in toric manifolds considered by Batyrev and Borisov
\cite{BB} as candidates for mirror symmetry. Each fibration in $\CC$
has a Lagrangian section.

An anti-symplectomorphism of a symplectic manifold $X$ with
symplectic form $\omega$ is a self-diffeomorphism $\phi$ of $X$ such
that
\[
\phi^*\omega = -\omega.
\]
An anti-symplectomorphism $\phi$ of $X$ such that $\phi^2 = \id_X$
is called an anti-symplectic involution.

In the present paper, we define a class of Lagrangian sections $\Cc$
for fibrations $f : X \rightarrow B$ of class $\CC.$ For each
fibration $f \in \CC,$ there exists at least one section $\sigma \in
\Cc.$ Our main results are the following.
\begin{thm}\label{tm:m}
Let $f: X \rightarrow B$ be a Lagrangian fibration of class $\CC.$
Let $\sigma$ be a Lagrangian section of $f$ of class $\Cc.$ There
exists a unique anti-symplectic involution $\phi_{f,\sigma}$ of $X$
such that
\begin{equation}\label{eq:2ct}
f\circ \phi_{f,\sigma} = f, \qquad \phi_{f,\sigma} \circ \sigma =
\sigma.
\end{equation}
That is, there exists a unique anti-symplectic involution
$\phi_{f,\sigma}$ of $X$ preserving the fibers of $f$ and fixing the
section $\sigma.$ Assuming existence, uniqueness continues to hold
for an arbitrary Lagrangian section $\sigma.$
\end{thm}
\begin{thm}\label{tm:tw}
Let $f : X \rightarrow B$ be a Lagrangian fibration of class $\CC$
and let $\sigma_0,\sigma_1,$ be two Lagrangian sections of class
$\Cc.$ There exists a unique symplectomorphism $t : X \rightarrow X$
satisfying
\begin{equation}\label{eq:tw}
f\circ t = f, \qquad t \circ \sigma_0 = \sigma_1.
\end{equation}
Assuming existence, uniqueness continues to hold for arbitrary
Lagrangian sections of $f.$
\end{thm}

There are several motivations for proving Theorem \ref{tm:m}.
Recently, there has been considerable research devoted to defining
Gromov-Witten type invariants for symplectic manifolds equipped with
an anti-symplectic involution \cite{W1,W2,S1}. Moreover, in the
presence of an anti-symplectic involution, an open-string mirror
correspondence was found \cite{RSW}. Most known examples of
anti-symplectic involutions come from real algebraic geometry.
Theorem~\ref{tm:m} constructs a vast number of examples of
symplectic manifolds with anti-symplectic involutions in a purely
symplectic way.  In Section~\ref{ssec:m}, we give a conjectural
mirror symmetry interpretation of Theorem~\ref{tm:m} that explains
conditions \eqref{eq:2ct}. In \cite{S2}, Theorem~\ref{tm:m} is
applied to show unobstructedness and calculate Lagrangian Floer
cohomology for smooth fibers of $f.$

Theorem~\ref{tm:tw} is important in the proof that
Theorem~\ref{tm:m} holds for any section $\sigma \in \Cc.$ Moreover,
as we explain in Section~\ref{ssec:m}, Theorem~\ref{tm:tw} has a
mirror symmetry interpretation of independent interest. In
Section~\ref{ssec:ns}, Theorem~\ref{tm:tw} allows us to construct
anti-symplectic involutions that do not fix a section of a
fibration.

\subsection{Mirror symmetry of symmetries}\label{ssec:m}
We briefly review some aspects of mirror symmetry necessary to put
Theorems~\ref{tm:m} and~\ref{tm:tw} in context. We discuss two
conjectures and the evidence in their favor.

\subsubsection{The Hodge diamond}
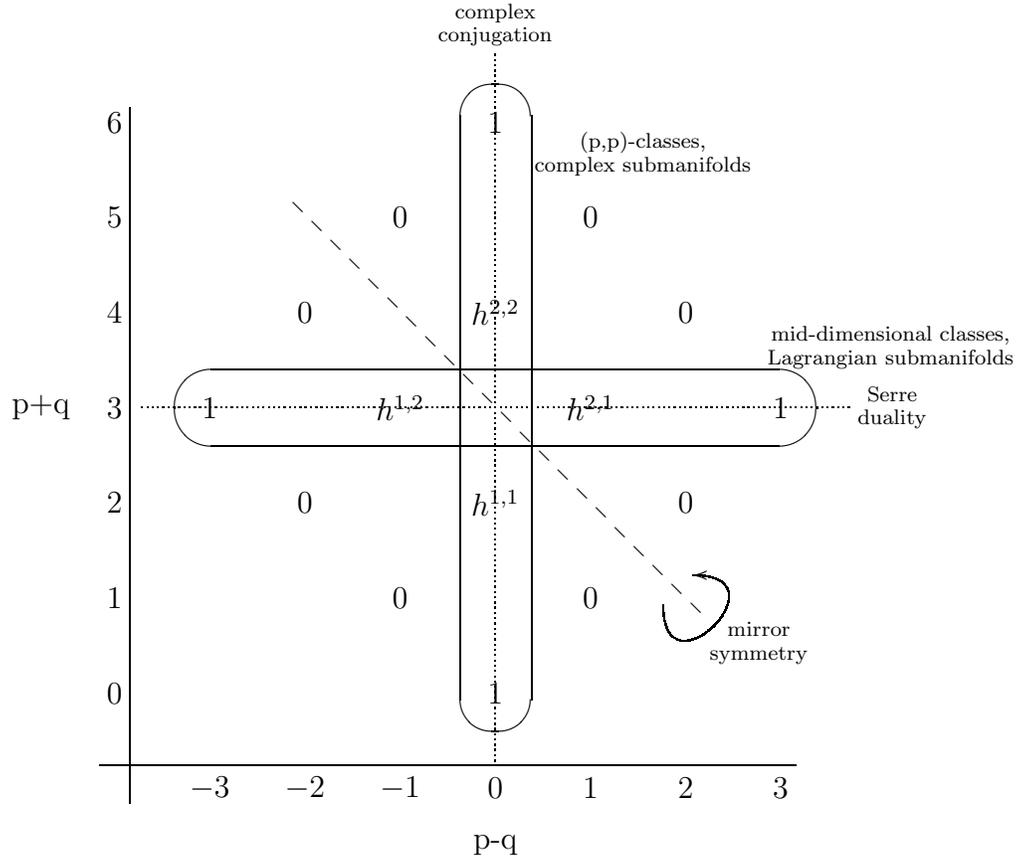
\begin{figure}
\centering
\[
\begin{xy}
{ \entrymodifiers={+++[o]}
 \xymatrix"D" @=3pc @!0{
6 & & & & 1 & & & & \\
5 & & & 0 & & 0 & & & \\
4 & & 0 & & h^{2,2} & & 0 \\
3 & 1 & & h^{1,2} & & h^{2,1} & & 1 \\
2 & & 0 & & h^{1,1} & & 0 \\
1 & & & 0 & & 0 & & & \\
0 & & & & 1 & & & & \\
{} & -3 & -2 & -1 & 0 & 1 & 2 & 3
 }}
,"D1,1"+<2mm,2mm>;"D8,1"+<2mm,-2mm>**@{-} %
,"D8,1"+<-2mm,3mm>;"D8,8"+<2mm,3mm>**@{-} %
,"D4,2"-<9mm,0mm>;"D4,8"+<9mm,0mm>*+!L\txt\tiny{Serre \\ duality}**@{.} %
,"D1,5"+<0mm,9mm>*+!D\txt\tiny{complex \\conjugation};"D7,5"-<0mm,9mm>**@{.} %
,"D2,3"+<-1.5mm,2mm>;"D6,7"+<2mm,-2mm>*+!UL\txt\tiny{mirror \\ symmetry}**@{--} %
,{"D6,7"+<-3mm,-1mm> \ar@(d,r)"D6,7"+<1mm,3mm>}%
,"D1,5";"D7,5"**\frm<44pt>{-}  %
,"D4,2";"D4,8"**\frm<44pt>{-} %
,"D1,5"+<4mm,0mm>*+!UL\txt\tiny{(p,p)-classes, \\complex
submanifolds} %
,"D4,8"+<-3mm,4mm>*+!DL\txt\tiny{mid-dimensional classes,
\\ Lagrangian submanifolds}
,"D4,1"+L*+!R\txt{p+q} %
,"D8,5"+D*+!U\txt{p-q}
\end{xy}
\]
\caption{The Hodge diamond} \label{fig:hd}
\end{figure}

A Calabi-Yau manifold is a K\"{a}hler manifold with trivial
canonical bundle. Mirror symmetry predicts that there exist pairs of
Calabi-Yau manifolds (X,Y) such that symplectic geometry on $Y$
mirrors complex geometry on $Y.$ Set $n = \dim_\C X = \dim_\C Y.$ A
concrete prediction of mirror symmetry is that
\begin{equation}\label{eq:ms}
H^q(X,\Omega^p_X) \simeq H^q(Y,\Omega^{n-p}_Y).
\end{equation}
Namely, the Hodge diamond of $Y$ is the reflection of the Hodge
diamond of $X$ about a diagonal. See Figure~\ref{fig:hd}, which
illustrates the case when $n = 3$ and $\pi_1(X) = \{1\}.$ We use the
notation $h^{p,q} = \dim H^q(X,\Omega^p_X).$

It follows from the isomorphism \eqref{eq:ms} that deformations of
the K\"{a}hler class of $X,$ which are classified by
$H^1(X,\Omega^1_X),$ are isomorphic to deformations of the complex
structure on $Y,$ which are classified by $H^1(Y,\Omega^{n-1}_Y).$
The middle-dimensional cohomology of $X,$ which contains the
Poincar\'{e} duals of Lagrangian submanifolds, is isomorphic to the
$(p,p)$ classes of $Y,$ which contain the Poincar\'{e} duals of
complex submanifolds.

Recall that the Hodge diamond of any K\"{a}hler manifold has two
symmetries: Serre duality and complex conjugation. Naturally, the
mirror isomorphism \eqref{eq:ms} preserves these symmetries of the
Hodge diamond. It is interesting to note, however, that mirror
symmetry exchanges complex conjugation and Serre duality. In the
remainder of this section, we trace the exchange of symmetries
through successively more refined descriptions of mirror symmetry.

\subsubsection{Homological mirror symmetry}

Homological mirror symmetry \cite{Kontsevich} can be seen as a
categorification of mirror symmetry on the level of Hodge diamonds.
Mirror symmetry on the level of Hodge diamonds implies an
isomorphism of vector spaces between the middle-dimensional
cohomology of $X$ and the $(p,p)$-classes of $Y:$
\[
\bigoplus_{p+q = n} H^{q}(X,\Omega_X^p) \ris \bigoplus_p
H^{p}(Y,\Omega_Y^p).
\]
Homological mirror symmetry replaces each vector space with a
category, and asserts an equivalence of categories.

Let $X$ be a symplectic manifold with $\dim_\R(X) = 2n.$ We would
like to replace $H^{n}(X)$ with a category. Since Lagrangian
submanifolds of $X$ have half the dimension of $X,$ it is natural to
look for a category with objects Lagrangian submanifolds. In fact,
from $X$ we can construct the $A^\infty$ category $Fuk(X)$
\cite{FOOO}. An object of $Fuk(X),$ is a Lagrangian submanifold $L
\subset X$ equipped with a unitary local system $E \rightarrow L,$ a
grading $\theta$ (see Section~\ref{sec:g}) and a $Pin$ structure
$\p.$ Depending on context, we may omit several of the data
comprising an object of $Fuk(X)$ from our notation when it does not
cause confusion. Morphisms between two objects $(L_1,E_1)$ and
$(L_2,E_2)$ are given by the Floer cohomology groups with local
coefficients $HF^*((L_1,E_1),(L_2,E_2)).$ From $Fuk(X),$ one can
construct a triangulated category $\DF(X),$ as explained in
\cite{Kontsevich}. In general, it seems necessary to enlarge
$\DF(X)$ further \cite{Kontsevich}. We denote the enlargement as
well by $\DF(X).$

On the other hand, let $Y$ be a K\"{a}hler manifold. To $Y,$ we can
associate the triangulated category $\DC(Y),$ the derived category
of coherent sheaves on $Y.$ Perhaps the simplest objects of
$\DC(Y),$ are the structure sheaves of complex submanifolds. The
Poincar\'e duals of complex submanifolds all belong to $\oplus_p
H^{p}(Y,\Omega^p).$ According to the Hodge conjecture, the
Poincar\'e duals of complex submanifolds should generate all
rational $(p,p)$-classes. Thus, it makes sense to replace the vector
space $\oplus_p H^{p}(Y,\Omega^p_Y)$ with the category $\DC(Y).$

Suppose $X$ and $Y$ are Calabi-Yau manifolds, i.e. Kahler manifolds
with trivial canonical bundle. Homological mirror symmetry
\cite{Kontsevich} predicts that for certain pairs $(X,Y)$ there
exists an equivalence of triangulated categories
\[
m : \DF(X) \ris \DC(Y).
\]
Such pairs are called mirror pairs and $Y$ is called a mirror of
$X.$

In homological mirror symmetry, symmetries of a vector space should
be replaced with auto-equivalences of a category. As for any smooth
algebraic variety, the functor
\begin{equation}\label{eq:D}
\D := \rsh(-,\O_Y) : \DC(Y)
\ris \DC(Y)^{op}
\end{equation}
induces an equivalence of categories, and
\begin{equation}\label{eq:minv}
\D^{op} \circ \D \simeq \id.
\end{equation}
The auto-equivalence $D$ is closely related to the Serre duality
symmetry of the Hodge diamond. It is natural to ask whether the
functor mirror to $\D,$
\[
\D^\vee : = (m^{op})^{-1} \circ \D \circ m : \DF(X) \ris
\DF(X)^{op},
\]
is quasi-isomorphic to a geometrically defined functor. One goal of
this paper is to construct a geometric functor
\[
\I : \DF(X) \ris \DF(X)^{op}
\]
for which it is reasonable to conjecture that
\[
\I \simeq \D.
\]
However, before proceeding further, we pause to summarize what is
conjectured about the geometry of $m.$

\subsubsection{The SYZ conjecture}

The SYZ conjecture \cite{SYZ} takes a first step toward giving a
geometric interpretation of the homological mirror symmetry functor
$m.$ For each point $y \in Y,$ let $\sky_y$ denote the skyscraper
sheaf at $y.$ According to \cite{SYZ}, the functor $m^{-1}$ should
carry $\sky_y$ to a Lagrangian torus $L_y \subset X$ equipped with a
flat unitary line bundle $E_y \rightarrow L_y.$ One motivation for
this conjecture was that $RHom(\sky_y,\sky_y) \simeq
\Lambda^*(T_yY),$ while it is reasonable to conjecture that
\[
HF^*((L_y,E_y),(L_y,E_y)) \simeq H^*(L_y) \simeq
\Lambda^*(H^1(L_y)).
\]
Indeed, the first isomorphism would follow if the spectral sequence
computing Floer cohomology \cite{FOOO} degenerates at the $E_2$
term.

Furthermore \cite{SYZ}, the family of tori $L_y$ should completely
fill $X$ and lead to a Lagrangian fibration $f:X \rightarrow B,$
where $B$ is a three-dimensional manifold. This fibration may have
singular fibers. Conversely, given a Lagrangian fibration $f:
X\rightarrow B,$ it should be possible to construct a mirror $Y_f,$
as the moduli space of pairs $(L,E),$ where $L$ is a fiber of $f$
and $E \rightarrow L$ is a flat unitary line bundle. The Lagrangian
fibrations constructed in \cite{CB-M} provide concrete examples in
which the SYZ conjecture can be tested.

The SYZ conjecture is sharpened by Fukaya in \cite{Fuk}. From now on
let $f:X \rightarrow B$ denote a Lagrangian fibration and let $Y_f$
be the corresponding mirror. Let $\sigma : B \rightarrow X$ be a
Lagrangian section of $f.$ It is suggested in \cite{Fuk} that the
choice of $f$ should play an important role in an intrinsic
construction of the equivalence of categories $\DF(X) \ris
\DC(Y_f).$ Indeed, given a coherent sheaf $\F$ on $Y_f$ and $y \in
Y_f,$ let $\F_y$ denote the fiber of $\F$ at $y$ and let $\sky_y$
denote the skyscraper sheaf at $y.$ We have canonically,
\begin{equation}\label{eq:st}
\F_y^\vee \simeq Hom(\F,\sky_y).
\end{equation}
Assume the object $(L,E)$ of $Fuk(X)$ is mirror to $\F,$ and let
$(L_y,E_y)$ be the fiber of $f$ mirror to $\sky_y.$ It follows from
homological mirror symmetry and isomorphism \eqref{eq:st} that we
must have an isomorphism
\[
\F_y^\vee \simeq HF^0((L,E),(L_y,E_y)).
\]
That is, we may calculate the fibers of the sheaf $\F$ on $Y_f$ only
knowing its mirror $L.$ A family version of Floer homology should
piece the fibers of $\F$ together to give $\F$ itself. In
\cite{Fuk2}, the choice of a Lagrangian section $\sigma$ plays an
important role in piecing together the fibers. The mirror functor
sends the Lagrangian $\sigma$ to the structure sheaf of $Y_f.$

In summary, the mirror functor depends on the choices of a
Lagrangian fibration $f$ and a Lagrangian section $\sigma.$ From now
on, we include $f$ and $\sigma$ in our notation for the mirror
functor
\[
m_{f,\sigma} : \DF(X) \ris \DC(Y_f).
\]
Similarly, we write
\[
\D^\vee_{f,\sigma} = (m_{f,\sigma}^{op})^{-1} \circ \D \circ
m_{f,\sigma}.
\]

\subsubsection{Duality Conjecture}

Suppose $f: X \rightarrow B$ is a Lagrangian fibration of class
$\CC.$ Let $\phi_{f,\sigma}$ be the involution of $X$ given by
Theorem~\ref{tm:m}. Let $(L,E,\theta,\p)$ be an object of $\DF(X).$
We define
\begin{equation}\label{eq:I}
\I_{f,\sigma}((L,E,\theta,\p)) =
(\phi_{f,\sigma}(L),\phi_{f,\sigma*}E^\vee,-\theta\circ\phi_{f,\sigma}^{-1},\phi_{f,\sigma*}\p).
\end{equation}
A short calculation shows that formula \eqref{eq:I} defines a
functor
\[
\I_{f,\sigma} : \DF(X) \longrightarrow \DF(X)^{op}.
\]
See \cite{Na}. The reversal of morphisms results from the fact that
$\phi_{f,\sigma}$ changes the sign of the symplectic form. For
signs, see \cite{S2}. The following conjecture has appeared in
various forms throughout the mirror-symmetry literature. See for
example Arinkin-Polishchuk in the case of the elliptic curve
\cite{AP} and Nadler in the case of the cotangent bundle \cite{Na}.
In both cases, the conjecture is a theorem.
\begin{con}\label{cj:}
The geometrically defined functor $\I_{f,\sigma}$ is
quasi-isomorphic to $\D^\vee_{f,\sigma}.$
\end{con}

We briefly present some evidence for Conjecture~\ref{cj:}. First of
all, it is clear from the definition that $\I_{f,\sigma}^{op} \circ
\I_{f,\sigma} \simeq \id$ by analogy with equation~\eqref{eq:minv}.
Furthermore, just as $\D$ preserves the structure sheaf of $Y_f,$ so
too $\I_{f,\sigma}$ preserves the Lagrangian section $\sigma$
equipped with the trivial local system. In Section~\ref{sec:g}, we
explain the appropriate choice of grading for $\sigma.$

A distinctive property of $\D$ is that it preserves $\sky_y$ up to a
shift in grading by $\dim_\C Y.$  It follows that a geometric
functor $\I_{f,\sigma}$ isomorphic to $D^\vee_{f,\sigma}$ should
preserve $(L_y,E_y)$ up to a shift in grading by $\dim_\C X =
\dim_\C Y.$ Indeed, Theorem~\ref{tm:m} guarantees that for a fiber
$L_y$ of $f,$ we have $\phi_{f,\sigma}(L_y) = L_y.$ Suppose $L_y$ is
smooth and hence a torus. By Corollary~\ref{cor:m},
$\phi_{f,\sigma}$ acts on $L_y$ by the inverse map of the torus
group. So, if $E_y$ is a flat unitary line bundle then
$\phi_{f,\sigma*}E_y = E_y^\vee.$ It follows that
$\phi_{f,\sigma*}E_y^\vee = E_y.$ In Section~\ref{sec:g} we verify
under some reasonable assumptions that for the natural choice of
grading $\theta_y$ on a torus fiber $L_y,$ $\I_{f,\sigma}$ shifts
$\theta_y$ by $\dim_\C X.$ Thus $\I_{f,\sigma}$ preserves the mirror
of $\sky_y$ up to a shift by $\dim_\C X = \dim_\C Y.$

Finally, we state a theorem concerning the derived category of
coherent sheaves that mirrors the uniqueness claim of
Theorem~\ref{tm:m}. The proof appears in Section~\ref{sec:cs}. In
the following, $\D$ is the functor defined in \eqref{eq:D}.

\begin{thm}\label{tm:o1}
Let $Y$ be a smooth projective variety of dimension $n.$ Let $\D' :
\DC(Y) \rightarrow \DC(Y)^{op}$ denote an equivalence of categories
such that
\[
\D'(\sky_y) \simeq \sky_y[n], \;\; \forall y \in Y, \qquad \qquad
\D'(\O_Y) = \O_Y.
\]
Then $\D' \simeq \D.$
\end{thm}

\subsubsection{Twist conjecture}
Let $\L \rightarrow Y_f$ denote a holomorphic line bundle over
$Y_f.$ Tensoring with $\L$ defines an auto-equivalence
\[
\T : \DC(Y) \ris \DC(Y).
\]
As before, we define the mirror auto-equivalence $\T^\vee$ by
\[
\T^\vee = (m_{f,\sigma})^{-1} \circ \T \circ m_{f,\sigma} : \DF(X)
\ris \DF(X).
\]
Let $(\sigma_\L,E_\L)$ denote the Lagrangian submanifold with
unitary local system mirror to $\L$ by the mirror isomorphism
$m_{f,\sigma}.$

We assume that $\sigma_\L$ is a Lagrangian section of $f:
X\rightarrow B$ and $E_\L$ has rank one. Furthermore, we assume that
$\sigma_\L \in \Cc.$ By Theorem~\ref{tm:tw}, there exists a unique
symplectomorphism $t : X \rightarrow X$ such that $f \circ t = f$
and $t \circ \sigma = \sigma_\L.$ Let $\widehat E_\L = f^*f_*E_\L,$
which is a flat unitary line bundle on $X$ since $\sigma_\L$ is a
section. Define an auto-equivalence $\t$ of $\DF(X)$ by
\[
\t((L,E,\theta,\p)) = (t(L),(t_*E)\otimes\widehat E_\L,\theta\circ
t^{-1},t_*\p).
\]
We make the following conjecture.
\begin{con}\label{con:tw}
The geometrically defined functor $\t$ is quasi-isomorphic to
$\T^\vee.$
\end{con}
Previously, Kontsevich described the functor $\t$ in terms of
monodromy transformations arising from complex structure moduli of
$X.$ See \cite{Hor,ST}. Given the right Lagrangian section
$\sigma_\L,$ the functor $\t$ can be used to construct the
homogeneous coordinate ring of $Y$ as described in \cite{Zas}.

We briefly present some evidence in favor of
Conjecture~\ref{con:tw}. Let $E_0 \rightarrow \sigma$ denote the
trivial rank-$1$ unitary local system. Just as
\[
\T(\O_Y) = \L,
\]
so too, since $t(\sigma) = \sigma_\L$ and $\widehat
E_\L|_{\sigma_\L} = E_\L,$ we have
\[
\t((\sigma,E_0)) = (\sigma_\L,E_\L).
\]
Also, just as
\[
\T(\sky_y) \simeq \sky_y, \qquad \forall y \in Y,
\]
so too, since $t(L_y) = L_y$ and $\widehat E_\L|_{L_y}$ is
trivial, we have
\[
\t((L_y,E_y)) \simeq (L_y, E_y).
\]
Finally, we state a theorem for the derived category of coherent
sheaves that mirrors the uniqueness claim of Theorem~\ref{tm:tw}.
The proof appears in Section~\ref{sec:cs}.
\begin{thm}\label{tm:o2}
Let $Y$ be a smooth projective variety. Let $\T'$ be an
auto-equivalence of $\DC(Y)$ such that
\[
\T'(\sky_y) \simeq \sky_y, \;\;\forall y \in Y, \qquad \qquad
\T'(\O_Y) \simeq \L.
\]
Then $\T' \simeq \T.$
\end{thm}

Finally, here is a justification of our assumption that $\sigma_\L$
is a section and $E_\L$ has rank $1.$ By homological mirror
symmetry, we have
\[
HF^*((\sigma_\L,E_\L),(L_y,E_y)) \simeq RHom(L,S_y) \simeq
\L_y^\vee.
\]
Since the fiber $\L_y$ is one dimensional, we conclude that
$HF^*((\sigma_\L,E_\L),(L_y,E_y))$ is one dimensional. The Floer
complex is a direct sum of tensor products of the fibers of the
local coefficient systems at intersection points of the Lagrangian
submanifolds. So, it is natural to assume that $\sigma$ intersects
each fiber $L_y$ at one point and $\rk E_\L = 1.$

\subsection{Idea of proof}
We briefly explain the idea of our proof of Theorems~\ref{tm:m}
and~\ref{tm:tw}. First suppose $f_0: X_0 \rightarrow B_0$ is a
Lagrangian fibration that is a smooth submersion. Let $\sigma_0$ be
a smooth Lagrangian section.

Recall that the cotangent bundle $T^*B_0$ has a canonical symplectic
form. Let $Z$ denote the zero section of $T^*B_0$. Let $\pi : T^*B_0
\rightarrow B_0$ denote the canonical projection. We also use $\pi$
to denote the induced projection of quotients of $T^*B_0.$

\begin{prop}[\cite{Dui}]\label{pr:aa}
There exists a unique lattice bundle $\Lambda_0 \subset T^*B_0$ and
a unique symplectomorphism
\[
\Theta: T^*B_0/\Lambda_0 \rightarrow X_0
\]
such that
\begin{equation}\label{eq:2cc}
\Theta \circ Z = \sigma_0, \qquad f_0 \circ \Theta   = \pi.
\end{equation}
\end{prop}
\begin{proof}[Sketch of proof]
Let $b \in B_0$ and let $\xi$ be a cotangent vector to $B_0$ at $b.$
Choose a function $h$ on $B_0$ such that $dh|_b = \xi.$ Define $H =
h\circ f_0.$ Let $\Phi_H$ be the time-one map of the Hamiltonian
flow of $H.$ Define $\widetilde \Theta : T^*B_0 \rightarrow X_0$ by
\[
\widetilde \Theta(b,\xi) = \Phi_H\circ\sigma_0(b).
\]
Define $\Lambda_0 = \widetilde \Theta^{-1}(\sigma_0).$ It is not
hard to check that $\widetilde \Theta$ descends to define a map
$\Theta$ on the quotient $T^*B_0/\Lambda_0$ with the required
properties.

To check uniqueness, assume $\Lambda_0' \subset T^*B_0$ is a lattice
bundle and
\[
\Theta' : T^*B_0/\Lambda_0' \rightarrow X_0
\]
is a symplectomorphism such that $\Theta' \circ Z = \sigma_0$ and
$f_0 \circ \Theta'   = \pi.$ Let $\widetilde \Theta'$ denote the
composition of $\Theta'$ with the quotient map $T^*B_0 \rightarrow
T^*B_0/\Lambda_0'.$ Let $b,\xi$ and $h,$ be as above. Define a
function $K$ on $T^*B_0$ by $K = h \circ \pi.$ Let $\Phi_K$ be the
time-one map of the Hamiltonian flow of $K.$ Since
\[
H \circ \widetilde \Theta' = h \circ f_0 \circ \widetilde \Theta' =
h \circ \pi = K,
\]
and $\widetilde \Theta'$ is a symplectomorphism, we conclude that
\begin{equation}\label{eq:HK}
\widetilde \Theta'\circ \Phi_K = \Phi_H \circ \widetilde \Theta'.
\end{equation}
An explicit calculation shows that
\begin{equation}\label{eq:bxi}
\Phi_K(b,0) = (b,\xi).
\end{equation}
Using equations \eqref{eq:HK} and \eqref{eq:bxi}, we conclude
\[
\widetilde \Theta'(b,\xi)  = \widetilde \Theta' \circ \Phi_K \circ Z(b) \\
 = \Phi_H \circ \widetilde \Theta' \circ Z(b) \\
 = \Phi_H \circ \sigma_0(b) \\
 = \widetilde \Theta(b,\xi).
\]
Thus $\widetilde \Theta' = \widetilde \Theta.$ It follows
immediately that $\Lambda_0' = \Lambda_0.$
\end{proof}

\begin{cor}\label{cor:m}
There exists a unique anti-symplectic involution $\phi_0$ of $X_0$
such that
\begin{equation}\label{eq:2c}
f_0 \circ \phi_0 = f_0, \qquad \phi_0 \circ \sigma_0 = \sigma_0.
\end{equation}
\end{cor}
\begin{proof}
Let $-\id$ denote the anti-symplectic involution of $T^*B_0$ given
by negative the identity transformation on fibers. We also use
$-\id$ to denote the induced involution of quotients of $T^*B_0.$
Let $\Theta$ denote the symplectomorphism constructed in
Proposition~\ref{pr:aa}. The diffeomorphism $\phi_0$ of $X_0$
defined by
\begin{equation}\label{eq:p0}
\phi_0 = \Theta \circ (-\id) \circ \Theta^{-1}.
\end{equation}
is an anti-symplectic involution satisfying conditions
\eqref{eq:2c}. We have proved existence.

To prove uniqueness, let $\phi_0'$ denote any anti-symplectic
involution satisfying conditions \eqref{eq:2c}. It follows that
\[
\phi_0' \circ \Theta \circ (-\id) : T^*B_0/\Lambda_0 \rightarrow X_0
\]
is a symplectomorphism satisfying conditions \eqref{eq:2cc}. By the
uniqueness claim of Proposition~\ref{pr:aa}, we conclude $\phi_0'
\circ \Theta \circ (-\id) = \Theta,$ and consequently $\phi_0' =
\phi_0.$
\end{proof}

\begin{cor}\label{cor:tw}
Let $\sigma_0$ and $\sigma_0'$ be two Lagrangian sections of $f_0.$
There exists a unique symplectomorphism $t_0 : X_0 \rightarrow X_0$
such that
\begin{equation}\label{eq:2ccc}
f_0 \circ t_0 = f_0, \qquad t_0 \circ \sigma_0 = \sigma_0'.
\end{equation}
\end{cor}
\begin{proof}
Let $\Theta$ (resp. $\Theta'$) denote the symplectomorphism given by
applying Proposition~\ref{pr:aa} to $\sigma_0$ (resp. $\sigma_0'$).
Clearly,
\[
t_0 =  \Theta'\circ \Theta^{-1} : X_0 \rightarrow X_0
\]
is a symplectomorphism satisfying conditions \eqref{eq:2ccc}.

To prove uniqueness, let $t_0' : X_0 \rightarrow X_0$ be any
symplectomorphism satisfying conditions \eqref{eq:2ccc}. Then
\[
t_0' \circ \Theta : T^*B_0/\Lambda_0 \rightarrow X_0
\]
satisfies conditions \eqref{eq:2cc} for section $\sigma'.$ It
follows that $t_0' \circ \Theta = \Theta',$ and consequently $t_0' =
t_0.$
\end{proof}

With this background, we outline the proofs of Theorems~\ref{tm:m}
and~\ref{tm:tw}. Suppose $f:X \rightarrow B$ is a Lagrangian
fibration of class $\CC.$ The construction of $f$ given in
\cite{CB-M} realizes $X$ as the compactification of an open dense
submanifold $X_0 \subset X$ such that $f_0 := f|_{X_0}$ is a
Lagrangian fibration that is a smooth submersion with a Lagrangian
section $\sigma_0.$ To obtain $X$ from $X_0,$ local models of
singular fibrations are glued onto $X_0$ matching up $\sigma_0$ to
sections of the local models. The main technical part of this paper
is devoted to constructing a fiber preserving anti-symplectic
involution fixing a section on each local model. By the denseness of
$X_0$ in $X,$ Corollary~\ref{cor:m} guarantees that all local
involutions piece together to form the global involution $\phi$ of
Theorem~\ref{tm:m}. Uniqueness of $\phi$ follows similarly. The same
approach together with Corollary~\ref{cor:tw} is used to construct
and prove the uniqueness of $t$ as in Theorem~\ref{tm:tw}.

\subsection{Involutions that do not fix a section}\label{ssec:ns}
A simple example of Calabi-Yau manifold with an anti-symplectic
involution is the Fermat quintic,
\[
Q = \{(z_0,\ldots,z_4) \in \C P^4 | \sum_{i = 0}^4 z_i^5 = 0\},
\]
with the anti-symplectic involution $\phi_Q$ induced by complex
conjugation. It is easy to see that the fixed locus of $\phi_Q$ is
connected. In fact, it is homeomorphic to $\R P^3.$

On the other hand, the fixed locus of any involution constructed by
Theorem~\ref{tm:m} cannot be connected. In particular, any
anti-symplectic involution fixing a section cannot have connected
fixed locus for the following reason. By the proof of
Corollary~\ref{cor:m}, if $\dim_\R X = 2n,$ the intersection of the
fixed locus with any smooth fiber is $2^n$ points. The section is
clearly a component of the fixed locus, but it intersects each fiber
in only one point, thus it cannot be the whole fixed locus.

In Remark~\ref{rem:c} below, we show how to construct
anti-symplectic involutions that do not fix a Lagrangian section.

The proof of the following proposition appears in
Section~\ref{sec:as}.
\begin{prop}\label{pr:ac}
Let $f: X \rightarrow B$ be of class $\CC.$ If $\phi_{f}: X
\rightarrow X$ is an anti-symplectomorphism satisfying $f \circ
\phi_f = f$ then
\begin{equation}\label{eq:sq}
\phi_f^2 = \id_X.
\end{equation}
In other words, all fiber-preserving anti-symplectomorphisms are
involutions. Furthermore, if $t: X \rightarrow X$ is a
symplectomorphism such that $f \circ t = f$ then the following
equation holds
\begin{equation}\label{eq:ac}
\phi_{f} \circ t^{-1} = t \circ \phi_{f}.
\end{equation}
\end{prop}
\begin{rem}\label{rem:1}
Let $f: X \rightarrow B$ be of class $\CC$ and $\sigma_0 \in \Cc.$
Let $\sigma_1$ be an arbitrary Lagrangian section such that there
exists a symplectomorphism $t: X \rightarrow X$ satisfying
conditions \eqref{eq:tw}. Then
\[
\phi_{f,\sigma_1} := t \circ \phi_{f,\sigma_0} \circ t^{-1}
\]
is an anti-symplectic involution satisfying conditions
\eqref{eq:2ct}. Moreover, by Proposition~\ref{pr:ac}, we have
\[
\phi_{f,\sigma_1} = t^2 \circ \phi_{f,\sigma_0}.
\]
\end{rem}

\begin{rem}\label{rem:c}
Let $f: X \rightarrow B$ be of class $\CC$ and $\sigma_0,\sigma_1
\in \Cc.$ Let $t$ be as in Theorem~\ref{tm:tw}. It follows from
Proposition~\ref{pr:ac} that $\phi' : = t \circ \phi_{f,\sigma_0}$
is an anti-symplectic involution. It follows from our definition of
$\Cc$ that if $\phi'$ fixes a Lagrangian section $\sigma_2$ then
$\sigma_2 \in \Cc.$ In this case, Remark~\ref{rem:1} applied to
$\sigma_0$ and $\sigma_2$ implies that $t$ has a square root.
Conversely, if for some reason $t$ does not have a square root, we
conclude that $\phi'$ does not fix a Lagrangian section.
\end{rem}
The following corollary is an immediate consequence of
Proposition~\ref{pr:ac} and Theorems~\ref{tm:m} and~\ref{tm:tw}.
\begin{cor}\label{cor:ns}
Let $f : X \rightarrow B$ be a Lagrangian fibration of class $\CC$
and let $\sigma_0,\sigma_1,$ be two Lagrangian sections of class
$\Cc.$ There exists a unique anti-symplectic involution $\phi : X
\rightarrow X$ satisfying
\begin{equation*}
f\circ \phi = f, \qquad \phi \circ \sigma_0 = \sigma_1.
\end{equation*}
Assuming existence, uniqueness continues to hold for arbitrary
Lagrangian sections of $f.$
\end{cor}
\begin{rem}
The proof in Section~\ref{sec:as} implies that
Proposition~\ref{pr:ac} continues to hold for any Lagrangian
fibration that is a smooth submersion with a smooth Lagrangian
section. Corollary~\ref{cor:ns} holds for any Lagrangian fibration
that is a smooth submersion and any smooth Lagrangian sections.
\end{rem}

\subsection{Acknowledgements}
We would like to thank P. Seidel for introducing us to
Conjecture~\ref{cj:} and for many helpful comments. We would like to
thank R. Bezrukavnikov for his patient help with the proof of
Corollary~\ref{cor:maa}. We would also like to thank K. Fukaya,
D. Kazhdan, and G. Tian for their helpful comments. R. Casta\~no-Bernard was
partially supported by a K-State FDA award 2291 and IHES. D. Matessi was partially supported by MIUR (``Riemannian metrics and differentiable manifolds'', PRIN 05). J. Solomon
was partially supported by NSF grant DMS-0703722.

\section{Lagrangian fibrations}\label{sec:bk}

\medskip

In \cite{CB-M}, the first two authors provide a method to construct
Lagrangian torus fibrations of $6$-dimensional symplectic manifolds homeomorphic to known Calabi-Yau manifolds. Recall that an integral affine structure $\mathscr A$
on a topological manifold is an atlas of charts whose change of coordinate maps
are affine maps with integral linear part, i.e. elements of $\R^n \rtimes \operatorname{SL}(\Z,n)$. The basic idea is to start with an integral affine manifold with singularities, $(B,\Delta,\mathscr A)$, where the $B$ is a topological $n$-manifold such that $B_0=B-\Delta$ has an integral affine structure $\mathscr A$. Here the discriminant $\Delta$ has codimension $2$ and the affine
structure is assumed to be \emph{simple}. Roughly speaking,
simplicity means that around points of $\Delta$, $B_0$ is locally affine isomorphic to given models of integral affine manifolds, satisfying certain natural
properties, such as having unipotent monodromy, cf. \cite{CB-M}
Definition 3.14 for details.

The affine structure on $B_0=B-\Delta$ induces a family of maximal
lattices $\Lambda\subseteq T^\ast B_0$, together with a symplectic
manifold $X_0$ and an exact sequence

\[
0\rightarrow\Lambda\rightarrow T^\ast B_0\rightarrow X_0\rightarrow
0.
\]

This gives us a Lagrangian $T^n$ bundle $f_0:X_0\rightarrow B_0$.
The manifold $X_0$ can be compactified to a topological $n$-manifold
$X$ by gluing on Gross' local models of topological $T^n$ fibrations
\cite{TMS}. To define a symplectic structure on $X$, in other words,
to achieve a \emph{symplectic} compactification of $X_0$, one needs
local models of \emph{Lagrangian} fibrations with singular fibres.
These models were studied in \cite{RCB1}, \cite{CB-M-torino},
\cite{CB-M-stitched}. In dimension $n=2$, $\Delta$ consists of a
finite collection of points and the symplectic compactification of
$X_0$ is achieved by gluing a standard model of a Lagrangian
fibration over a disc with a nodal central fibre; this model is
known in symplectic geometry as a simple focus-focus fibration. This
construction gives compact symplectic 4-manifolds with Lagrangian
2-torus fibrations. For a specific choice of integral affine $S^2$
with 24 singularities, one obtains a symplectic 4-manifold
diffeomorphic to a K3 surface (cf. \cite{LeungSym} or \cite{CB-M} Theorem 3.22). In
dimension $n=3$, $\Delta$ is typically a graph with trivalent
vertices, labeled either positive or negative. In this case the affine
structure around edges or positive and negative vertices is isomorphic
to the one induced on the base of models of three different kinds of local Lagrangian fibrations: respectively the so-called \textit{generic}, \textit{positive} and \textit{negative} fibrations (see Section~\ref{sec:local}).
The models for generic and positive fibrations can be regarded as 3-dimensional analogues of focus-focus fibrations; in particular, the fibrations are $T^2$-invariant, have
codimension 2 discriminant and are given by smooth fibration maps.
On the other hand, the model for a negative fibration is
$S^1$-invariant, the fibration map is piecewise smooth and its
discriminant locus has mixed codimension 1 and 2. This model can be
regarded as a perturbation Gross' topological version of the
negative fibration used in \cite{TMS}. This perturbation forces the
discriminant locus to drop codimension in a small neighborhood of a
negative vertex. As a consequence, the compactification of
$f:X_0\rightarrow B_0$ achieved in \cite{CB-M},  required a
modification of $\Delta$, more precisely, a fattening near the
negative vertices.

\subsection{The class $\mathcal{C}$.}\label{sec:C}

Given a simple integral affine $3$-manifold with singularities $(B,
\Delta, \mathscr A)$ a \textit{localized thickening} of $\Delta$ is
given by the data $(\Delta^{\blacklozenge}, \{ D_{p^-} \}_{p^- \in
\mathcal N})$ where:
\begin{itemize}
\item[(i)] $\Delta^{\blacklozenge}$ is the closed subset obtained from $\Delta$ after replacing a
           neighborhood of each negative vertex with a  localized amoeba, i.e. an amoeba as in Figure \ref{fig: amoeba} after the end of each leg is pinched down to dimension one.
\item[(ii)] $\mathcal N$ is the set of negative vertices and for each $p^-\in \mathcal N$, $D_{p^-}$ is a $2$-disk containing the codimension $1$ component of $\Delta^{\blacklozenge}$ around  $p^-$.
\end{itemize}
Given a localized thickening, define
\[
B_{\blacklozenge} = B - \left( \Delta \cup \bigcup_{p^- \in \mathcal
N} D_{p^-} \right),
\]
and denote by $\mathscr A_{\blacklozenge}$ the restriction of the
affine structure on $B_{\blacklozenge}$. Let $X_{\blacklozenge} = T^{\ast}B_{\blacklozenge} / \Lambda$ with standard symplectic form and $f_{\blacklozenge}: X_{\blacklozenge} \rightarrow B_{\blacklozenge}$ be the projection.

\medskip
The main result of \cite{CB-M} is the following:

\medskip
\begin{thm}\label{class_C}
Given a compact simple integral affine $3$-manifold with
singularities $(B, \Delta, \mathscr A)$, satisfying some additional
mild hypothesis, there is a localized thickening
$(\Delta_{\blacklozenge}, \{ D_{p^-} \}_{p^- \in \mathcal N})$  and
a smooth, compact symplectic $6$-manifold $X$ together with a
piecewise smooth Lagrangian fibration $f: X \rightarrow B$ such that
\begin{itemize}
\item[(i)] $f$ is smooth except along $\bigcup_{p^- \in \mathcal N} \,  f^{-1}(D_{p^{-}})$;
\item[(ii)] the discriminant locus of $f$ is $\Delta_{\blacklozenge}$;
\item[(iii)] there is a commuting diagram
\begin{equation*}
\begin{CD}
X_{\blacklozenge}  @>\Psi>> X\\
@Vf_\blacklozenge VV  @VVfV\\
B_{\blacklozenge} @>\iota>> B
\end{CD}
\end{equation*}
where $\psi$ is a symplectomorphism and $\iota$ the inclusion;
\item[(iv)] over a neighborhood of a positive vertex of $\Delta_{\blacklozenge}$ the fibration is
             positive, over a neighborhood of a point on an edge the fibration is generic-singular, over
             a neighborhood of $D_{p^-}$ the fibration is Lagrangian negative.
\item[(v)] $f$ has a section, $\sigma$, such that $\sigma (B)\subset X$ is a smooth Lagrangian submanifold such that $\sigma (B)\cap\Crit f=\varnothing$.
\end{itemize}
\end{thm}

Compact symplectic manifolds $X$ and Lagrangian fibrations
$f:X\rightarrow B$ as in Theorem \ref{class_C} define a class,
$\mathcal C$. The class $\mathcal C$ includes a large number of
symplectic models of mirror pairs of Calabi-Yau manifolds with SYZ
dual Lagrangian fibrations. Examples of these include the quintic
3-fold and its mirror and Batyrev-Borisov pairs of Calabi-Yau
manifolds (cf.\cite{CB-M} for details).

\medskip
In the next Section, we review the local models used for the
compactification in Theorem \ref{class_C}, and provide a
case-by-case proof of existence of anti-symplectic involutions for
each model.

\section{Local existence of involutions}\label{sec:local}

\subsection{Focus-focus fibrations}
We show how to construct fibre-preserving anti-symplectic
involutions in dimension $n=2$. Consider the case of proper
fibrations with focus-focus type singularities.  Here is a simple
example:

\begin{ex} \label{smooth:ff}
Let $X = \C^{2} - \{ z_1 z_2 + 1 = 0 \}$ and let $\omega$ be the
restriction to $X$ of the standard symplectic form on $\C^2$. The
following map $f:X\rightarrow\R^2$ is a Lagrangian fibration:
\begin{equation} \label{sm:ff}
f(z_1,z_2)=\left( \frac{|z_1|^2-|z_2|^2}{2}, \, \log |z_1z_2 +1|
\right).
\end{equation}
The only singular fibre is $f^{-1}(0)$, which is nodal with one node
at the point $(0,0)$. Clearly conjugation $\phi:(z_1, z_2) \mapsto
(\bar z_1, \bar z_2)$ on $\C^2$ is a fibre-preserving
anti-symplectic involution. The fixed locus of $\phi$ is $\R^2 - \{
x_1x_2 + 1 = 0 \}$ which has $3$ connected components. Two of them
(i.e. the connected components of $\{ x_1x_2 +1 < 0 \}$) are
sections of $f$ not containing the singular point. The other one
(i.e. the set $\{ x_1x_2 +1 > 0 \}$) is mapped $2$ to $1$ by $f$
except at $(0,0)$ which is a branched point.
\end{ex}

We now describe the construction of a general nodal fibration.
Details can be found in \cite{RCB}. First, let us discuss a local model for the singularity. The standard focus-focus singularity is the (non-proper) map $q: \C^2 \rightarrow
\C$ given by
\begin{equation} \label{local_model}
 q(z_1, z_2) = z_1 \bar z_2.
\end{equation}
Here $z_1 = y_1+ i y_2$, $z_2 = x_1 + i x_2$ and the symplectic form
is $\omega = \sum dx_j \wedge dy_j$. The real and imaginary parts of
the map $q$ are $q_1 = x_1y_1 + x_2y_2$ and $q_2 = x_1y_2 - x_2y_1$
respectively. If $v_{q_j}$ denotes the Hamiltonian vector field
corresponding to $q_j$ and $g_j^t$ its flow, we have that
\[ g_1^t(z_1, z_2) = (e^{-t}z_1, e^t z_2), \]
\[ g_2^t(z_1, z_2) = (e^{it} z_1, e^{it}z_2). \]
Notice that $v_{q_2}$ induces an $S^1$ action. In fact,  if $\tau =
e^{- t_1 + it_2}$ the we have
\begin{equation}  \label{ff:action1}
 (g_1^{t_1} \circ g_2^{t_2}) (z_1, z_2) = ( \tau \, z_1,  \bar \tau^{-1} \, z_2)
\end{equation}
which gives a $\C^{\ast}$ action.

If $B  = \{ b=b_1 + i b_2, \, |b| < 1 \}$, we restrict the above map
$q$ to $Y = q^{-1}(B)$. We have two Lagrangian sections of $q$,
$\Sigma_j: B \rightarrow Y$, $j=1,2$, given by $\Sigma_1(b) = (1,
\bar b)$ and $\Sigma_2(b) = (b, 1)$. Define the maps
\begin{eqnarray*}
 \phi_j: \C^{\ast} \times B & \rightarrow & Y \\
      (\tau, b) & \mapsto & \tau \cdot \Sigma_{j}(b)
\end{eqnarray*}
Clearly $\phi_j$ describes the orbit of the section $\Sigma_j$ via
the Hamiltonian flow. Now let
\begin{equation} \label{v1}
 V_1 = \{ (\tau,  b) \in \C^{\ast} \times B \, | \, |b| < |\tau| < 1 \}
\end{equation}
\begin{equation} \label{v2}
 V_2 = \{ ( \tau, b) \in \C^{\ast} \times B \, | \, 1 < |\tau| < |b|^{-1} \},
\end{equation}
and let $U_j = \phi_j(V_j)$. Denote $U = (U_1 \cup U_2)$. Clearly $U
\cup \{ (0,0) \}$ is an open neighborhood of the singular point
$(0,0)$.

Now, suppose we are given a proper nodal fibration $f: X \rightarrow
B$ with singular point $p \in X$, and a Lagrangian section $\sigma$
of $f$ with $\sigma(B)\subset X^{\#} = X - p$. We can describe
$X^\#$ in terms of action angle-coordinates using a non-proper
version of Proposition \ref{pr:aa}. Consider  $B \subset \C$ the
unit disk and $T^{\ast}B$ with its standard symplectic form. For any
smooth function $H: B \rightarrow \R$ define the following closed
1-forms on $B$:
\begin{eqnarray*}
\lambda_1 & = & - \log |b| \ db_1 + \Arg b \ db_2 + dH \\
\lambda_2 & = & 2 \pi  \, db_2,
\end{eqnarray*}
and consider the integral lattice $\Lambda \subset T^{\ast}B$
spanned by $\lambda_1$ and $\lambda_2$ and let $J_H = T^{\ast}B /
\Lambda$. Then one can prove (cf. \cite{RCB1}{Thm. 2.5 and 3.1})
that with a suitable choice of $H$, there is a unique
fibre-preserving symplectomorphism
\begin{equation}\label{th:non-prop}
\Theta: J_H\rightarrow X^\#
\end{equation}
which maps the zero section to $\sigma$. Uniqueness of $\Theta$ follows from Proposition \ref{pr:aa} and by continuity.

Conversely, we now show how to use the above descriptions to construct a proper nodal fibration. Any such fibration can be obtained by gluing the neighborhood $U \cup \{ (0,0) \}$ of the focus-focus singularity to the space $J_H$ defined above. First of all,
notice that $J_H$ also has a $\C^{\ast}$ action, namely for every
$\alpha \in T^{\ast}_bB$,
\begin{equation} \label{ff:action2}
 \tau \cdot (b, \alpha) = (b, \, -\log |\tau| db_1 + \Arg \tau \, db_2 + \alpha)
\end{equation}
which are just translations along the fibres. Now let $L_1$ be the
Lagrangian section in $J_H$ given by the graph of $dH$ and let
$L_2$ be the one given by the zero section. We can define maps
\begin{eqnarray*}
 \psi_j: \C^{\ast} \times B & \rightarrow & J_H \\
      (\tau, b) & \mapsto & \tau \cdot L_{j}(b)
\end{eqnarray*}

Notice that $ \psi_1(1,b) = L_1(b)$ and $\psi_1(b,b) = L_2(b)$. Take
$V_j \subseteq \C^{\ast} \times B$, $j=1,2$ as in (\ref{v1}) and
(\ref{v2}) and denote $U_j^{\prime} = \psi_j(V_j)$ and $U^{\prime}=
U_1^{\prime} \cup U_2^{\prime}$. Now define $g:
U^{\prime} \rightarrow U$ by $g|_{U_j^{\prime}} = \phi_j
\circ \psi^{-1}_{j}$. We can use $g$ to glue the focus-focus singularity to $J_H$,
to form a symplectic manifold with a nodal fibration.

All nodal fibrations can be realized this way for any given function
$H$ on $B$. Furthermore, given two nodal fibrations of $X$ and $X'$
determined by functions $H$ and $H'$ there is a fibre preserving
symplectomorphism in a neighborhood of the singular fibres if an
only if $H$ and $H'$ have the same germ at the origin. In other
words, the invariant of a germ of nodal fibration is a germ of a
function on $B$ at the origin, a formal power series in two
variables. Moreover, this power series is independent of the choice
of Lagrangian section of $f$ (cf. \cite{RCB1}).

\medskip
We have the following useful result:

\begin{lem}\label{lem:tw:non-proper}
Let $\sigma_1$ and $\sigma_2$ be two Lagrangian sections of a nodal fibration $f:X\rightarrow B$ not intersecting the singular point $p\in X$. There exists a unique symplectomorphism $t : X \rightarrow X$
such that
\begin{equation}\label{eq:2ccc-nonprop}
f \circ t = f, \qquad t \circ \sigma_1 = \sigma_2.
\end{equation}
\end{lem}
\begin{proof}
For the section $\sigma_1$ there is a function $H$ on $B$ and a map
$\Theta:J_{H}\rightarrow X^\#$ which maps the zero section to
$\sigma_1$ (see above).  The section $\sigma_2$ corresponds to a
section $\sigma'$ in $J_H$. Now we can define a symplectomorphism
$t'$ on $J_H$ which is given by translation by $\sigma'$ on the
fibres (clearly $t'$ maps the zero section to $\sigma'$). With this,
one can show that $t^\# = \Theta \circ t' \circ \Theta^{-1}$ extends
to a symplectomorphism $t$ of $X$ satisfying the required
conditions. In fact, if $\sigma'(b) = s_1(b) db_1 + s_2(b) db_2$,
then one can describe $t'$ in terms of the $\C^*$ action by
\[ t'(b,\alpha) = \tau(b) \cdot (b, \alpha), \]
where $\tau(b) = e^{- s_1(b) + i s_2(b)}$ (cf.. (\ref{ff:action2})).
Since the map $g$, which glues the singularity to $J_H$, also
matches the $\C^{\ast}$ actions, we have that $t'$ corresponds to
the following map on the local model for the focus-focus singularity
(see (\ref{ff:action1})):
\[ t'(z_1, z_2) = (\tau(b) z_1, \bar{\tau}(b)^{-1} z_2) \]
where $b = q(z_1, z_2) = z_1 \bar z_2$. Clearly this map extends
smoothly to the singularity. To prove uniqueness of $t$, one
restricts $t$ to $X_0$ and applies Corollary \ref{cor:tw}.
\end{proof}

Finally we have:

\begin{prop}\label{invo:nodal}
Let $f:X \rightarrow B$ be a nodal fibration with a Lagrangian section $\sigma$ not intersecting the singular point $p\in X$. There exists a unique anti-symplectic involution $\iota_{f,\sigma}$ of $X$ such that
\begin{equation}\label{eq:nodal-inv}
f \circ \iota_{f,\sigma} = f, \qquad \iota_{f,\sigma} \circ \sigma = \sigma.
\end{equation}
\end{prop}

\begin{proof}
Consider first the local model $q: \C^{2} \rightarrow \C$ given by
(\ref{local_model}). Clearly the map $\iota: (z_1, z_2) \mapsto
(\bar z_2, \bar z_1)$ is a fibre-preserving anti-symplectic
involution, which exchanges the Lagrangian sections $\Sigma_1$ and
$\Sigma_2$. Similarly, $J_H$ has the anti-symplectic involution
\[ \iota: (b, \alpha) \mapsto (b, dH(b) - \alpha), \]
which exchanges $L_1$ and $L_2$ and fixes $\frac{1}{2}L_1$.  Since the gluing map $g$ satisfies
$g \circ \iota = \iota \circ g$, $\iota$ extends to a
fibre-preserving anti-symplectic involution on $(X,\omega, f)$ fixing the section $\sigma':=\frac{1}{2}L_1$.

Now suppose that $\sigma$ is another section of $f$ not intersecting the singular point.
>From Lemma~\ref{lem:tw:non-proper} we know that there is a fibre-preserving symplectomorphism $t$ sending $\sigma$ to $\sigma'$. Then the map $\iota_{f,\sigma} = t^{-1} \circ \iota \circ t$ is a fiber preserving anti-symplectomorphism with the required properties. To prove uniqueness one simply restricts $f$ and $\sigma$ to $X_0$ and $B_0$, respectively, and applies Corollary \ref{cor:m}.
\end{proof}

\subsection{Generic-singular fibration.} \label{generic:singular} An almost identical construction can be carried out for generic-singular fibrations in
dimension $6$. The (non-proper) local model for the singularity is
the map $q: \C^{2} \times S^1 \times (0,1) \rightarrow \C \times
S^1$ given by
\[ q: (z_1, z_2, e^{i \theta}, r) \mapsto (z_1 \bar z_2, r), \]
which is singular along $\Crit{q} = \{ 0 \} \times S^1 \times
(0,1)$. The Hamiltonian flow induces a $\C^{\ast} \times S^1$ action
on $\C^2 \times S^1$ given by
\[ (\tau, e^{is}) \cdot (z_1, z_2, e^{i \theta}, r)= (\tau z_1, \bar \tau^{-1} z_2, e^{i (\theta + s)}, r). \]

The space $X^{\#}$ is constructed as follows. Let $D \subset \C$ be
the open unit disk and $B = D \times (0,1)$. Given a smooth function
$H$ on $B$, we form the lattice $\Lambda \subseteq T^{\ast}B$
generated by the periods
\begin{eqnarray*}
\lambda_1 & = & - \log |b| \ db_1 + \arg b \ db_2 + dH \\
\lambda_2 & = & 2 \pi  \, db_2, \\
\lambda_2 & = & dr.
\end{eqnarray*}
Then $X^{\#} = T^{\ast}B / \Lambda$ also has a $\C^{\ast} \times
S^1$ action
\[ (\tau, e^{is}) \cdot (b, \alpha) =
                 (b, \, -\log |\tau| db_1 + \Arg \tau \, db_2 + s dr + \alpha). \]
We can now glue a neighborhood of $\Crit(q)$ to $X^{\#}$ using this
action just like in the $4$ dimensional case. This gives a
symplectic manifolds $X$ and a proper Lagrangian fibration over a
cylinder $B=D\times (0,1)$ with discriminant $\Delta=\{0\}\times
(0,1)$ and singular fibres being a product of $S^1$ and a nodal
fibre. This construction provides all generic-singular fibrations.
All generic singular fibrations can be realized using the above
method and, as in the nodal case, germs of generic-singular
fibrations are classified \cite{RCB1}. The proof of the following
two statements is the same as the proof of Lemma
\ref{lem:tw:non-proper} and Proposition \ref{invo:nodal}:

\begin{lem}\label{lem:tw:non-proper-gen}
Let $\sigma_1$ and $\sigma_2$ be two Lagrangian sections of a generic-singular fibration $f:X\rightarrow B$ not intersecting the singular set $\Crit f \subset X$. There exists a unique symplectomorphism $t : X \rightarrow X$
such that
\begin{equation}\label{eq:2ccc-nonprop-gen}
f \circ t = f, \qquad t \circ \sigma_1 = \sigma_2.
\end{equation}
\end{lem}

\begin{prop}\label{invo:gen}
Let $f:X \rightarrow B$ be a generic-singular fibration with a Lagrangian section $\sigma$ not intersecting the critical locus $\Crit f\subset X$. There exists a unique anti-symplectic involution $\iota_{f,\sigma}$ of $X$ such that
\begin{equation}\label{eq:gen-inv}
f \circ \iota_{f,\sigma} = f, \qquad \iota_{f,\sigma} \circ \sigma = \sigma.
\end{equation}
\end{prop}

In this case the fixed point locus consists of seven components, six
of which are sections. We leave the details of the proof to the
reader.

\subsection{Positive fibration} The situation is completely analogous to the case of fibrations of nodal type discussed above. First we give
an explicit example.
\begin{ex} \label{proper:pos}
Let $X=\C^3-\{ 1+ z_1z_2z_3=0\}$ with standard symplectic structure.
Let $f:X\rightarrow\R^3$ be given by
$$f(z_1,z_2,z_3) = (\log |1+z_1z_2z_3|, |z_1|^2-|z_2|^2, |z_1|^2-|z_3|^2).$$
The map $f$ defines a proper Lagrangian fibration having the
topology of a positive fibration. The singular fibres lie over a
trivalent vertex:
\[
\Delta=\{ b_1=0, b_2=b_3\geq 0\}\cup\{ b_1=b_2=0,b_3\leq 0\}\cup\{
b_1=b_3=0, b_2\leq 0\}.
\]
The fibres over the edges have generic-singular type discussed above
while the central fibre is homeomorphic to $S^1\times T^2$ after a 2
cycle, $\{x\}\times T^2$, is collapsed to $x\in S^1$. It is easy to
see that $f$ is invariant under the anti-symplectic involution of
$X$ given by conjugation $\iota: (z_1, z_2, z_3) \mapsto (\bar z_1,
\bar z_2, \bar z_3)$. The fixed set of $\iota$ is the set $\R^3 - \{
1 + x_1x_2x_3 \}$, which has $5$ connected components. The one
containing $(0,0,0)$ is mapped by $f$ generically $4$ to $1$, while
the other four components are sections not intersecting the singular
locus.
\end{ex}

A neighborhood of the singular locus in the above example is modeled
on the following example due to Harvey and Lawson \cite{HL}.

\begin{ex}\label{ex. HL}
On $\C^3$ define
\begin{equation}
F(z_1,z_2,z_3) = (\im z_1z_2z_3, |z_1|^2-|z_2|^2, |z_1|^2-|z_3|^2)
\end{equation}
Here the last two components define the moment map $\mu$ of a
Hamiltonian $T^2$-action. The critical locus of $F$ is $\Crit
(F)=\bigcup_{ij}\{ z_i=z_j=0\}$ and its discriminant locus is
$\Delta$ as in Example~\ref{proper:pos}. The regular fibres are
homeomorphic to $\R\times T^2$. The singular fibre over $0\in\Delta$
is homeomorphic to $\R\times T^2$ after $\{ p\}\times T^2$ is
collapsed to $p\in\R$. All the other singular fibres are
homeomorphic to $\R\times T^2$ after a two cycle $\{ p \} \times T^2
\subset \R \times T^2$ is collapsed to a circle. The map $\iota:
(z_1, z_2, z_3) \mapsto (- \bar z_1, \bar z_2, \bar z_3)$ defines a
fibre preserving anti-symplectic involution. Notice that the smooth
part of every singular fibre has two connected components and
$\iota$ sends one to the other.
\end{ex}

In the above example, the Hamiltonian flow associated to the
components of $F$ induces an $\R \times T^2$ action on $\C^3$ which
is free and transitive on the smooth fibres. Let us denote by
$t=(t_1, t_2, t_3)$ the coordinates on $\R \times T^2$, where $t_2$
and $t_3$ are periodic (of period $1$) and by $(t, z) \mapsto t
\cdot z$ the action on $\C^3$.  Consider a neighborhood $B$ of
$(0,0,0) \in \R^3$ and let $Y = F^{-1}(B)$. We can find sections of
$F$, $\Sigma_j: B \rightarrow Y$, $j=1,2$, chosen so that $\iota
\circ \Sigma_1 = \Sigma_2$. There is a function $\tau: B \rightarrow
\R \times T^2$ such that $\tau(b) \cdot \Sigma_1 (b) = \Sigma_2(b)$.
Actually, $\tau$ is defined only on $B - \Delta$, and it is shown in
\cite{RCB1} that the first component $\tau_1$ of $\tau$ tends to
$+\infty$ as $b$ approaches the discriminant locus $\Delta$. Define
the maps
\begin{eqnarray*}
 \phi_j: (\R \times T^{2}) \times B & \rightarrow & Y \\
      (t, b) & \mapsto & t \cdot \Sigma_{j}(b).
\end{eqnarray*}
Notice that $\phi_1(0,b) = \Sigma_1(b)$ and $\phi_1(\tau(b), b) =
\Sigma_2(b)$ if $b \notin \Delta$. Since $\tau_1$ is big near
$\Delta$, we may assume $\tau_1(b) > 0$ for all $b \in B$. For
$j=1,2$ define subsets
\begin{eqnarray}
 V_1 & = & \{ (t,b) \in  (\R \times T^{2}) \times B \, | \, t_1 \in (0, \tau_1(b))\} \label{pos:v1}\\
 V_2 & = &\{ (t,b) \in  (\R \times T^{2}) \times B \, | \, t_1 \in (- \tau_1(b), 0) \}\label{pos:v2}
\end{eqnarray}
and let $U_j = \phi_j(V_j)$. We have that $U_1 \cup U_2$ is a $T^2$
invariant open set, whose closure $\bar U$ is a (closed)
neighborhood of $\Crit(F)$ having as boundary the $T^2$ orbit of
$\Sigma_1(B) \cup \Sigma_2(B)$.

Now, for any smooth function $H$ on $B$, define one forms on $B$ as
follows
\begin{align} \label{pos:per}
\lambda_1 = \sum_{j=1}^3 \tau_j \, db_j +dH ,\quad \lambda_2 =
db_2, \quad  \lambda_3 = db_3.
\end{align}
It can be shown that these are all closed one forms. Let $X^{\#} =
T^{\ast}B / \Lambda$, where $\Lambda$ is the integral lattice
generated by the $\lambda_j$'s. Also on $X^{\#}$ there is an $\R
\times T^2$ action given by
\[ t \cdot (b, \alpha) = (b, \, \sum_{j=1}^{3} t_j db_j + \alpha). \]
Now let $L_1$ be the Lagrangian section in $X^{\#}$ given by the
graph of $dH$ and let $L_2$ be the zero section. We can define maps
\begin{eqnarray*}
 \psi_j: (\R \times T^2) \times B & \rightarrow & X^{\#} \\
      (t, b) & \mapsto & t \cdot L_{j}(b)
\end{eqnarray*}
Notice that $ \psi_1(0,b) = L_1(b)$ and $\psi_1(\tau(b),b) =
L_2(b)$. Take $V_j$, $j=1,2$, as in (\ref{pos:v1}) and
(\ref{pos:v2}) and denote $U_j^{\prime} = \psi_j(V_j)$ and
$U^{\prime}= U_1^{\prime} \cup U_2^{\prime}$. We can now define a
map $g: U^{\prime} \rightarrow U$ such that $g|_{U_j^{\prime}} =
\phi_j \circ \psi^{-1}_{j}$ and use it to glue the singularity to
$X^{\#}$, to form $X$ with the positive fibration $f: X \rightarrow
B$. Also in this case we have:

\begin{lem}\label{lem:tw:non-proper-pos}
Let $\sigma_1$ and $\sigma_2$ be two Lagrangian sections of a positive fibration $f:X\rightarrow B$ not intersecting the singular set $\Crit f\subset X$. There exists a unique symplectomorphism $t : X \rightarrow X$
such that
\begin{equation}\label{eq:2ccc-nonprop-pos}
f \circ t = f, \qquad t \circ \sigma_1 = \sigma_2.
\end{equation}
\end{lem}

\begin{prop}\label{invo:pos}
Let $f:X \rightarrow B$ be a positive fibration with a Lagrangian section $\sigma$ not intersecting the singular set $\Crit f\subset X$. There exists a unique anti-symplectic involution $\iota_{f,\sigma}$ of $X$ such that
\begin{equation}\label{eq:gen-pos}
f \circ \iota_{f,\sigma} = f, \qquad \iota_{f,\sigma} \circ \sigma = \sigma.
\end{equation}
\end{prop}

The proofs are almost word by word like in Lemma~\ref{lem:tw:non-proper} and Proposition~\ref{invo:nodal}. As in Example~\ref{ex. HL}, the fixed point locus of the involution consists of $5$ connected components, four of which are sections.

\subsection{A piecewise smooth fibration}\label{sec:neg_pws}
Here we prove the existence of a fibre-preserving anti-symplectic
involution on interesting examples of piecewise smooth fibrations
used in \cite{CB-M}. We now recall the construction of this
fibration. Consider the following $S^{1}$ action on $\C^3$:
\begin{equation} \label{action}
e^{i\theta}(z_1,z_2, z_3)
          =(e^{i\theta}z_1,e^{-i\theta}z_2,z_3).
\end{equation}
It is Hamiltonian with moment map:
\begin{equation}\label{eq mu}
\mu (z_1, z_2, z_3)=\frac{|z_1|^2-|z_2|^2}{2}.
\end{equation}
The only critical value of $\mu$ is $t=0$. Now let $\gamma: \C^2
\rightarrow \C$ be the following $S^{1}$-invariant, piecewise smooth
map
\begin{equation}\label{eq. g}
\gamma (z_1,z_2)= \begin{cases}
\frac{z_1z_2}{|z_1|},\quad\text{when}\ \mu(z_1, z_2)\geq 0\\
\\
\frac{z_1z_2}{|z_2|},\quad\text{when}\ \mu(z_1, z_2) <0.
\end{cases}
\end{equation}
Define $\pi: \C^3 \rightarrow \C^2$ to be
\[ \pi(z_1, z_2, z_3) = (\gamma(z_1, z_2) , z_3) \]
and $\Log: (\C^{\ast})^2 \rightarrow \R^2$ to be
$$\Log(u_1, u_2) = (\log |u_1|, \log |u_2|).$$
It was shown in \cite{CB-M-torino} that given a symplectomorphism
$\Phi: \C^2 \rightarrow \C^2$ the map

\begin{equation}\label{s1inv:f}
 f = ( \mu, \Log \circ \Phi \circ \pi)
\end{equation}
defines a piece-wise smooth Lagrangian fibration on the open subset
of $\C^3$ given by
$$X = (\Phi \circ \pi)^{-1}((\C^{\ast})^{2}).$$

\medskip

Now consider the anti-symplectic involution $\iota: \C^3 \rightarrow
\C^3$ given by conjugation. It is easy to see that if $\Phi$
commutes with $\iota$ then $\iota(X) = X$ and $f \circ \iota = f$,
and therefore $\iota$ is fibre-preserving.

\begin{ex} \label{ex amoebous fibr}
If in the above construction we take
\begin{equation}\label{eq. linear sympl}
\Phi (u_1,u_2)=\frac{1}{\sqrt{2}}\left (
u_1-u_2,u_1+u_2-\sqrt{2}\right )
\end{equation}
then the map $f$ becomes
\begin{equation} \label{eq. the fibration}
f(z_1, z_2, z_3) = \left(\frac{1}{2}\left(|z_1|^2-|z_2|^2\right),
\log \frac{1}{\sqrt{2}}\left| \gamma - z_3\right|,
                  \log \frac{1}{\sqrt{2}} \left|\gamma + z_3 - \sqrt{2}\right|\right),
\end{equation}
The discriminant locus $\Delta$ is described as follows. Consider,
inside $(\C^\ast)^2$, the surface
\[ \Sigma = \{ v_1 + v_2 + 1 = 0 \}, \]
which is, topologically, a pair of pants. Then
\[ \Delta = \{ 0 \} \times \Log ( \Sigma). \]
Clearly $\Delta$ has the shape of an amoeba:

\begin{figure}[!ht]
    \centering
    \includegraphics[width=4cm,height=4cm,bb=0 0 258 258]{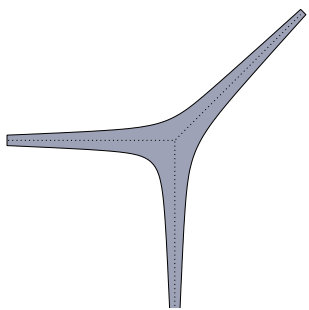}
    \caption{Amoeba of $v_1+v_2+1=0$}
    \label{fig: amoeba}
\end{figure}

For a discussion of the topology of the fibres in this example we
refer to \cite{CB-M-torino}. Observe that $\Phi$ commutes with the
conjugation map $\iota$, and therefore $\iota$ is fibre-preserving.
The fixed locus $S$ of $\iota$ is the complement in $\R^3$ of the
set
\[ K = \{ \gamma(x_1,x_2) - x_3 = 0 \} \cup \{ \gamma(x_1, x_2) + x_3 - \sqrt{2} = 0 \} \]
The reader may verify that $S = \R^3 - K$ has five connected
components, $S_1, S_2, \ldots, S_5$ containing respectively
$(0,0,2)$, $(0,0,-1)$, $(-1,-1, 1/2)$, $(1,1, 1/2)$, $(0,0,1/2)$.
Then $S_1$, $S_2$ and $S_5$ are mapped generically $2$ to $1$ while
$S_3$ and $S_4$ are sections.
\end{ex}

Now we verify that the same involution $\iota$ as above is
fibre-preserving also with respect to the version of the above
example where the legs are pinched down to codimension $2$ toward
the ends. Here is how we construct it.

\begin{ex}\label{ex amoebous_th fibr}
Consider the smooth function:
\[ H_0 = \frac{\pi}{4} \im (u_1 \overline{u}_2) \]
and let $\Phi_{H_0}$ be the Hamiltonian symplectomorphism associated
to $H_0$, i.e.
\[ \Phi_{H_0}: (u_1, u_2) \mapsto \frac{1}{\sqrt{2}} (u_1 - u_2, u_1 + u_2). \]
We now want a symplectomorphism which acts like $\Phi_{H_0}$ in a
small ball centered at the origin and like the identity outside a
slightly bigger ball. So choose a cut-off function $k: \R_{\geq 0}
\rightarrow [0,1]$ such that, for some $\epsilon > 0$,
\begin{equation} \label{gl:fn}
    k(t) = \left \{ \begin{array}{l}
                   1 \ \quad\text{when} \ 0< t \leq \epsilon; \\
                   0 \ \quad\text{when} \ t \geq 2\epsilon
                  \end{array} \right.
\end{equation}
 and define the Hamiltonian
\[ H = k(|u_1|^2 + |u_2|^2) H_0. \]
The Hamiltonian symplectomorphism $\Phi_H$ associated to $H$
satisfies
\[
\Phi_H(u_1, u_2 ) =
\begin{cases}
 \id_{\C^2},&\quad\text{when}\ |u_1|^2 + |u_2|^2 \geq 2\epsilon;\\
\\
\frac{1}{\sqrt{2}} (u_1 - u_2, u_1 + u_2),&\quad\text{when}\
                      |u_1|^2 + |u_2|^2 \leq \epsilon .
\end{cases}
\]
Now let $\Psi$ be the affine symplectomorphism
 \[
\Psi: (v_1,v_2) \mapsto
\frac{1}{\sqrt{2}}(v_1-v_2,v_1+v_2-\sqrt{2}).
\]
and finally, define $\Phi = \Psi \circ \Phi_H$. It is clear that

\begin{equation}\label{eq: Phi_thin_leg}
\Phi(u_1, u_2 ) =
\begin{cases}
 \Psi,&\quad\text{when}\ |u_1|^2 + |u_2|^2 \geq 2\epsilon;\\
\\
  ( - u_2, u_1 - 1),&\quad\text{when}\
                      |u_1|^2 + |u_2|^2 \leq \epsilon .
\end{cases}
\end{equation}

The fibration $f$ defined by (\ref{s1inv:f})  with this choice of
$\Phi$ has as discriminant locus a $3$-legged amoeba with the
horizontal leg pinched down to a line towards its ends.

Observe that $H$ satisfies $ H (\bar u_1, \bar u_2)= - H (u_1,
u_2)$. The reader may check that this property implies that $\Phi_H$
commutes with conjugation $\iota$. It follows that $\Phi$ commutes
with $\iota$ (since also $\Psi$ does).  Similarly one shows that the
symplectomorphism which pinches down all three legs at once is:

\begin{equation} \label{thin:fi}
\Phi(u_1, u_2 ) =
\begin{cases}
  ( - u_2, u_1 - 1),&\quad\text{when}\
                      |u_1|^2 + |u_2|^2 \leq \epsilon; \\
  \\
  (u_1 - 1, u_2 - \sqrt{2}), &\quad\text{when}\
                      |u_1|^2 + |u_2 - \sqrt{2}|^2 \leq \epsilon; \\
  \\
  \frac{1}{\sqrt{2}} ( u_1-u_2, \, u_1+u_2 ), &\quad\text{when}\
                       |u_2|^2 \geq M; \\
  \\
  \Psi,&\quad\text{everywhere else}.
\end{cases}
\end{equation}
which also commutes with $\iota$. So $\iota$ is fibre-preserving
also with respect to $f$ constructed with this $\Phi$.
\end{ex}

The fibration above has the right topology to be a good candidate
for a Lagrangian negative fibration. However, it fails to be smooth
over a large hyperplane containing the discriminant locus, and
therefore it is not a suitable model for the compactification as in
Theorem \ref{class_C}. A suitable model should be smooth away from a
disc containing the codimension 1 part of the discriminant. The
smoothing process is a delicate issue, and involves studying a class
of piece-wise smooth Lagrangian fibrations, called stitched
fibrations \cite{CB-M-stitched}.

\begin{rem}\label{rem:section-neg}
One may use lifts of sections of the $\Log$ fibration to obtain many Lagrangian sections of the fibration in Example \ref{ex amoebous fibr} not intersecting its singular locus (details cf. \cite{CB-M}{Proposition 5.9}). From all these, only the sections $S_3$ and $S_4$ are fixed by the involution given by complex conjugation. The same situation holds for the thin-legged Example \ref{ex amoebous_th fibr}. In what follows, we shall assume that our choice of a Lagrangian section $\sigma$ of Example \ref{ex amoebous_th fibr} is either given by $S_3$ or $S_4$.
\end{rem}

\subsection{Involutions of stitched fibrations}
This Section is rather technical, closely related to
\cite{CB-M-stitched}, and may be skipped on a first reading. The
main (new) results are Theorems \ref{thm:stitched:inv} and
\ref{thm:stitched_2:inv} where we provide conditions for the
existence of a smooth, fibre-preserving, anti-symplectic involution
of stitched fibrations, which will be used to prove the existence of
involutions of the negative fibration.

\begin{defi}\label{defi stitched}
Let $X$ be a smooth $2n$-dimensional symplectic manifold. Suppose
there is a free Hamiltonian $S^1$ action on $X$ with moment map
$\mu: X \rightarrow \R$. Let $X^+ = \{ \mu \geq 0 \}$ and $X^- = \{
\mu \leq 0 \}$. Given a smooth $(n-1)$-dimensional manifold $M$, a
map $f: X \rightarrow \R \times M$ is said to be a \textit{stitched
Lagrangian fibration} if there is a continuous $S^1$ invariant
function $G: X \rightarrow M$, such that the following holds:
\begin{itemize}
\item[\textit{(i)}] Let $G^{\pm} = G|_{X^{\pm}}$. Then $G^+$ and $G^-$ are restrictions of $C^\infty$ maps on $X$;
\item[\textit{(ii)}] $f$ can be written as $f =  (\mu, G)$ and $f$ restricted to $X^{\pm}$ is a proper submersion with connected Lagrangian fibres.
\end{itemize}
We call $Z = \mu^{-1}(0)$ the \textit{seam} and $\Gamma=f(Z)
\subseteq \{0\}\times M$ the \textit{wall}. We denote $f^\pm
=f|_{X^\pm}$.
\end{defi}


Denote $B = f(X)$ and $B^{\pm} = f(X^{\pm})$. In general, a stitched
fibration will only be piecewise $C^\infty$, however all its fibres
are smooth Lagrangian tori. Throughout this section we will always
assume (unless otherwise stated) that the pair $(B, \Gamma)$ is
diffeomorphic to the pair $(D^n, D^{n-1})$, where $D^k \subset \R^k$
is an open unit ball centered at the origin and $\R^{n-1}$ is
embedded in $\R^n$.

\medskip

Observe that the fibration in Example \ref{ex amoebous fibr}, when
restricted to $X-f^{-1}(\Delta)$, defines a stitched Lagrangian
fibration. In fact this is the main example. The seam is $Z =
\mu^{-1}(0) - f^{-1}(\Delta)$, notice that in this case $Z$ has
three connected components.

\medskip

The seam of a stitched fibration is an $S^1$-bundle $p:Z\rightarrow
\bar Z:=Z\slash S^1$ such that $f$ factors through $p$, i.e. we have
the diagram:
\[
\xymatrix{
Z \ar[dr]_{f|_Z} \ar[rr]^{p} &  & \bar{Z} \ar[dl]^{\bar{f}} \\
&  \Gamma }
\]
where $\bar{Z}$ has the reduced symplectic form and $\bar f$ is the
reduced Lagrangian fibration over the wall $\Gamma$.  We also have
the vertical $(n-1)$-plane distribution:
\[
\mathfrak{L}=\ker \bar f_\ast\subset T\bar Z
\]
tangent to the fibres of $\bar f$. In what follows we will define
certain invariants of the stitched fibration consisting of sections
of $\mathfrak{L}^*$ which are \textit{fiberwise closed}, in the
sense that they restrict to closed $1$-forms on the fibres of $\bar
f$.

\medskip
On the base of a stitched fibration we allow a more general set of
coordinates than just the smooth ones, which we define bellow.

\begin{defi}\label{defi:admissible}
A set of coordinates on $B \subseteq \R \times M$, given by a map
$\phi: B \rightarrow \R^n$, is said to be \textit{admissible} if the
components of $\phi=(\phi_1,\ldots ,\phi_n)$ satisfy the following
properties:
\begin{itemize}
\item[\textit{(i)}] $\phi_1$ is the restriction to $B$ of the projection map $\R \times M\rightarrow\R$;
\item[\textit{(ii)}]  for $j = 2, \ldots, n$ the restrictions of $\phi_j$ to $B^+$ and $B^-$ are locally
restrictions of smooth functions on $B$.
\end{itemize}
\end{defi}

Essentially, admissible coordinates are those such that $\phi \circ
f$ is again stitched. Let $f: X \rightarrow B$ be a stitched
Lagrangian fibration and let $\phi$ be a set of admissible
coordinates. For $j = 2, \ldots , n$,   $f_{j}^{\pm} = \phi_j  \circ
f|_{X^{\pm}}$ is the restriction of a $C^\infty$ function on $X$ to
$X^\pm$ and we can write $f=(\mu, f_2^\pm,\ldots ,f_n^\pm)$.

\medskip
Now we want to put stitched fibrations in a normal form. In the
smooth case, a proper Lagrangian submersion locally always admits
action-angle coordinates, defined up to the choice of a basis of
$H_1(F_b,\Z)$, where $F_b$ is the fibre over $b \in B$. In the case
of stitched fibrations we can generalize this idea as follows.
Assuming $B$ contractible, we choose a pair of bases
$\gamma^\pm=(\gamma_1,\gamma_2^\pm, \ldots ,\gamma_n^\pm)$ of
$H_1(X,\Z)$ such that

\begin{itemize}
\item[(a)] $\gamma_1$ is represented by an orbit of the $S^1$ action,
\item[(b)] $\gamma_j^{+} = \gamma_j^{-} + m_j \gamma_1$, for some $m_2, \ldots,m_n\in\Z$.
\end{itemize}

Condition (b) simply means that $p_\ast\gamma^+=p_\ast\gamma^-$
under the map $p_\ast:H_1(X,\Z)\rightarrow H_1(X/S^1,\Z)$. The
following proposition generalizes the notion of action angle
coordinates on the base.

\begin{prop}\label{prop:stitched_action} Let $f: X \rightarrow B$ be
a stitched fibration and let $\gamma^\pm$ be bases of $H_1(X,\Z)$
satisfying the above conditions. Then the restrictions of
$\gamma^\pm$ to $H_1(X^\pm, \Z)$ induce embeddings,
\[
\Lambda^\pm\hookrightarrow T^\ast_{B^\pm}.
\]
Let $\alpha^\pm: B^\pm\rightarrow\R^n$ be the corresponding action
coordinates satisfying $\alpha^\pm(b)=0$ for some $b \in \Gamma$.
Then the map
\[
\alpha=
\begin{cases}
\alpha^+ &\textrm{on}\ B^+\\
\alpha^-&\textrm{on}\ B^-
\end{cases}
\]
is an admissible change of coordinates. If $b_1,\ldots b_n$ denote
the action coordinates on $B$ given by $\alpha$, then $\{db_1,\ldots
db_n\}$ is a basis of $\Lambda^+$ and $\Lambda^-$.
\end{prop}

Recall that to establish the existence of action-angle coordinates,
in the classical case, one chooses a smooth Lagrangian section. In
the stitched case we choose a continuous section $\sigma: B
\rightarrow X$ such that $\sigma|_{B^{\pm}}$ are the restrictions of
smooth maps and $\sigma(B)$ is a smooth Lagrangian submanifold. Such
sections always exist locally, for instance in Example~\ref{ex
amoebous fibr} a component of the fixed locus of the anti-symplectic
involution is a section of this type. We denote a stitched fibration
$f: X \rightarrow B$ together with a choice of basis $\gamma$ of
$H_1(X, \Z)$ and a section $\sigma$ as above by $\mathcal F=(X, B,
f, \gamma, \sigma)$. When $W \subseteq B$ is an open set we usually
denote by $\mathcal F|_{W}$ the fibration $(f^{-1}(W), W, f, \gamma,
\sigma|_{W})$.

\begin{defi}\label{def:symp-eq}
Two stitched fibrations $(X, B, f, \gamma, \sigma)$ and $(X', B',
f', \gamma', \sigma')$, with seams $Z$ and $Z'$ respectively are
\textit{symplectically conjugate} if there are neighborhoods
$W\subseteq B$ of $\Gamma :=f(Z)$ and $W'\subseteq B'$ of $\Gamma':=
f'(Z')$, an $S^1$ equivariant $C^\infty$ symplectomorphism $\psi: f^{-1}(W) \rightarrow f'^{-1}(W')$
sending $Z$ to $Z'$ and a $C^\infty$ diffeomorphism $\phi: W \rightarrow W'$ such that:
$f' \circ \psi = \phi \circ f$, $\psi\circ\sigma=\sigma' \circ\phi$ and $\psi_\ast \gamma = \gamma'$.
The set of equivalence classes under this relation will be called
\textit{germs of stitched fibrations}.
\end{defi}

Notice that in the above definition we are allowed to shrink to a
smaller neighborhood of $\Gamma$ but not to a smaller $\Gamma$. So
germs are meant to be defined around $\Gamma$ and not around a
point.

\medskip
The following is a basic construction of stitched fibrations.

\begin{ex}[Normal forms] \label{ex:nf}
Let $(b_1,\ldots, b_n)$ be the standard coordinates on $\R^n$. Let
$(U,\Gamma)$ be a pair of subsets of $\R^n$ diffeomorphic to
$(D^n,D^{n-1})$ and $\Gamma = U\cap\{ b_1=0\}$. Define $U^+=U
\cap\{b_1\geq 0\}$ and $U^-=U\cap\{ b_1\leq 0\}$. Consider the
lattice $\Lambda= \spn\langle db_1,\ldots ,db_n\rangle_\Z$ and form
the symplectic manifold $T^\ast U\slash\Lambda$. Denote by $\pi$ the
standard projection onto $U$ and let $Z =\pi^{-1}(\Gamma)$. We may
consider the $S^1$ action on $T^\ast U\slash\Lambda$ given by
translations by multiples of $db_1$ in the fibres of $\pi$, whose
moment map is $\mu = b_1$. Suppose there is an open neighborhood
$V\subseteq T^\ast U\slash\Lambda$ of $Z$ and a map
$u:V\rightarrow\R^n$ which is a proper, smooth, $S^1$-invariant
Lagrangian submersion with components $(u_1,\ldots ,u_n)$ such that
$u|_{Z}=\pi$ and $u_1=b_1$. Now define the following subsets of
$T^\ast U\slash\Lambda$,
\[
Y^+ := \pi^{-1}(U^+),\quad Y:= Y^+ \cup V,\quad Y^- := Y \cap
\pi^{-1}(U^-)
\]
and define the map $f_u: Y \rightarrow \R^n$ by
\begin{equation} \label{u:st}
       f_u = \begin{cases}
                u \quad\text{on} \  Y^-, \\
               \pi \quad\text{on} \ Y^+.
       \end{cases}
\end{equation}
Clearly $f_u:Y\rightarrow \R^n$ is a stitched fibration. Denote
$B_u:=f_u(Y)$. The zero section $\sigma_0$ of $\pi$ is, perhaps
after a change of coordinates in the base, a section of $f_u$. Let
$\gamma_0$ be the basis of $H_1(Y,\Z)$ induced by $\Lambda$. We call
the stitched fibration $\mathcal F_u=(Y,B_u,f_u,\sigma_0,\gamma_0 )$
a \textit{normal form}.
\end{ex}
Now consider a normal form $\mathcal F_u$ and let $(b,y) = (b_1,
\ldots, b_n, y_1, \ldots, y_n)$ be canonical coordinates on
$T^{\ast}B_u$  so that $y$ gives coordinates on the fibre
$T_b^{\ast}B_u$. Let $W$ be a neighborhood of $\Gamma$ inside
$u(V)$. If $r \in \R$ is a parameter,  for any $b = (0, b_2, \ldots,
b_n) \in \Gamma$, let $(r,b)$ denote the point $(r,b_2, \ldots, b_n)
\in \R^n$. Given $(r,b) \in W$, denote by $L_{r,b}$ the fibre
$u^{-1}((r,b))$. For every fibre $F_b \subset Z$ of $\pi$, consider
the symplectomorphism
\begin{equation} \label{fb:symp}
 (y_1, \ldots, y_n, \sum_{k=1}^{n} x_k dy_k)
            \mapsto (x_1, b_2 + x_2, \ldots, b_n + x_n, y_1, \ldots, y_n),
\end{equation}
between a neighborhood of the zero section of $T^{\ast}F_b$ and a
neighborhood of $F_b$ in $V$. If $W$ is sufficiently small, for
every $(r,b) \in W$, the Lagrangian submanifold $L_{r,b}$ will be
the image of the graph of a closed $1$-form on $F_b$. Due to the
$S^1$ invariance of $u$ and the fact that $u_1=b_1$, this 1-form has
to be of the type
\[ r dy_1 + \ell(r,b), \]
where $\ell(r,b)$ is $S^1$ invariant, i.e. it may be considered as a
$1$-form on $\bar F_b := F_b / S^1$.  Denote by $\ell(r)$ the smooth
one parameter family of sections of $\mathfrak{L}^{\ast}$ such that
$\ell(r)|_{\bar{F}_b} = \ell(r,b)$. The condition $u|_Z=\pi$ implies
that $\ell (0,b)=0$. Furthermore, the $N$-th order Taylor series
expansion of $\ell(r)$ in the parameter $r$ can be written as
\begin{equation} \label{l:tay}
 \ell(r) =\sum_{k=1}^{N} \ell_k \, r^k + o(r^N),
\end{equation}
where the $\ell_k$'s are fiberwise closed sections of
$\mathfrak{L}^{\ast}$.

\begin{defi}\label{def:seq}
With the above notation, we define
\begin{itemize}
\item[(i)] $\mathscr L_{Z}$ the set of sequences
$\ell = \{\ell_k\}_{k\in\N}$ such that $\ell_k$ is a fiberwise
closed section of $\mathfrak L^\ast$;
\item[ii)] $\mathscr U_{Z}$ the set of pairs $(V,u)$ where $V\subseteq T^\ast U\slash\Lambda$ is a neighborhood of $Z$ and $u:V\rightarrow\R^n$ is a proper, smooth,  $S^1$-invariant Lagrangian submersion with components $(u_1,\ldots ,u_n)$ such that $u|_{Z}=\pi$ and $u_1=b_1$.
\end{itemize}
\end{defi}

As above, to a given $(V,u) \in \mathscr U_{Z}$ we can associate a
unique sequence $\ell \in \mathscr L_{Z}$. Conversely, in
\cite{CB-M-stitched}{\S 5} it is shown that for any given sequence
$\ell\in\mathscr L_{Z}$ there is some $(V,u)\in\mathscr U_{Z}$,
therefore a normal form, associated to it. Clearly, this $(V,u)$ is
not unique.

\medskip
In \cite{CB-M-stitched} the following result is proved:
\begin{prop}\label{prop:normalform} Every stitched fibration
$\mathcal F = (X, B, f, \sigma, \gamma)$ is symplectically conjugate
to a normal form $\mathcal F_u=(Y, B_u, f_u, \sigma_0, \gamma_0 )$
\end{prop}

When $\mathcal F$ is smooth, its normal form is $\mathcal F_\pi$,
this is Arnold-Liouville theorem. Given a stitched Lagrangian
fibration $\mathcal F=(X,B,f, \sigma, \gamma)$ with normal form
$\mathcal F_u=(Y, B_u, f_u, \sigma_0, \gamma_0 )$, we respectively
denote by $\znor$ and $\gnor$ the seam and the wall of $\mathcal
F_u$ and by $\zbnor$ the $S^1$ reduction of $\znor$.

\begin{defi}\label{def:seq-inv}
Let $\mathcal F=(X,B,f, \sigma, \gamma)$ be a stitched fibration
with normal form $\mathcal F_u=(Y, B_u, f_u, \sigma_0, \gamma_0 )$.
Let $\ell\in\mathscr{L}_{\zbnor}$ be the unique sequence determined
by $(V,u) \in \mathscr U_{\znor}$ defining $\mathcal F_u$. We call
$\inv (\mathcal F):=(\zbnor,\ell)$ the \textit{invariants} of
$\mathcal F$. We say that the invariants of $\mathcal F$ vanish if
for all $k\in\N$, $\ell_k\equiv 0$ when restricted to the reduced
fibres of $\mathcal F_u$. We say that the invariants of $\mathcal F$
are fiberwise constant if all the $\ell_k$'s are fiberwise constant.
\end{defi}
One can prove that $\inv (\mathcal F)$ is independent on the choice
of normal form. Moreover, we also have the following classification
results from \cite{CB-M-stitched}:

\begin{thm} \label{broken:constr2}
Given any pair $(U,\Gamma_\mathrm{nor})$ of subsets of $\R^n$,
diffeomorphic to $(D^n, D^{n-1})$ and with $\gnor = U \cap \{ b_1 =
0 \}$, a sequence $\ell = \{ \ell_k \}_{k \in \N} \in \mathscr
L_{\zbnor}$ and integers $m_2, \ldots, m_n$ such that
\begin{equation} \label{int:cond2}
  \int_{[db_j]} \ell_1 = m_j, \quad \text{for all}\ j=2, \ldots, n,
\end{equation}
there exists a smooth symplectic manifold $(X, \omega)$ and a
stitched Lagrangian fibration $f: X \rightarrow U$ satisfying the
following properties:
\begin{itemize}
\item[(i)] the coordinates $(b_1, \ldots, b_n)$ on $U$ are action coordinates of $f$ with $\mu = f^{\ast}b_1$ the moment map of the $S^1$ action;
\item[(ii)] the periods $\{ db_1, \ldots, db_n \}$, restricted to $U^{\pm}$ correspond to bases $\gamma^{\pm} = \{ \gamma_1, \gamma_2^{\pm}, \ldots, \gamma_n^{\pm} \}$ of $H_1(X, \Z)$ satisfying conditions (a) and (b) prior to Proposition \ref{prop:stitched_action};
\item[(iii)]  there is a Lagrangian section $\sigma$ of $f$, such that $(\zbnor, \ell)$ are the invariants of $(X,f, U, \sigma, \gamma^+)$.
\end{itemize}
\end{thm}

\begin{thm}\label{thm: grosso} Let $\mathcal F$ and $\mathcal F'$ be stitched fibrations. Then,
\begin{itemize}
\item[(i)] two stitched fibrations $\mathcal F$ and $\mathcal F'$ are symplectically conjugate
if and only if $\inv (\mathcal F)=\inv (\mathcal F')$;
\item[(ii)] $\mathcal F$ is smooth if and only if $\inv (\mathcal F)$ vanish;
\item[(iii)] $\mathcal F$ becomes smooth after an admissible change of coordinates on the base if and only if $\inv (\mathcal F)$ are fiberwise constant.
\end{itemize}
\end{thm}

In the above, fiberwise constant means that in the normal form, the
forms $\ell_k$ are independent of the $y$ coordinates. The set of
germs of stitched fibrations is therefore classified by the pairs
$(\zbnor, \ell)$. We say that a fibration is \textit{fake stitched}
if it becomes smooth after an admissible change of coordinates on
the base. The important consequence of Theorem~\ref{broken:constr2},
which was exploited in \cite{CB-M}, is that from a given set of
invariants we can form another one for example by summing to the
sequence $\ell$ another sequence or by multiplying elements $\ell_k$
by pull backs of smooth functions on the base. The new invariants
give rise to new stitched fibrations.

\medskip

\begin{rem}\label{rem:no-trans-stitched}
Observe that Corollary \ref{cor:tw} does not hold for stitched
fibrations in general, it only holds for stitched fibrations that
are fake. Indeed, let $f$ be a stitched fibration with a section
$\sigma$ and invariants $\ell_k$. Suppose $\sigma'$ is another
Lagrangian section of $f$ and assume there exists a
symplectomorphism $t$ of $X$ such that $f\circ t=f$ and
$t\circ\sigma=\sigma'$. This implies that $ \bar t^* \ell_k =\ell_k$
where $\bar t$ is the translation induced on the reduced fibration.
Therefore each $\ell_k$ is fiberwise constant, hence $f$ is smooth
after a suitable change of coordinates in the base.
\end{rem}

On a smooth Lagrangian fibration $f: X \rightarrow B$, with $B$
diffeomorphic to $D^n$, with a Lagrangian section $\sigma: B
\rightarrow X$, there always exists a unique fibre-preserving
anti-symplectic involution $\iota: X \rightarrow X$ fixing $\sigma$.
In fact, if $(b,y)$ are action-angle
coordinates on $X$ then we must have
\begin{equation} \label{iota:can}
\iota(b,y) = (b, -y)
\end{equation}

How about stitched fibrations? Do they admit smooth fibre-preserving
anti-symplectic involutions? First observe that given
$(X,B,f,\sigma, \gamma)$ with seam $Z$ and wall $\Gamma$, then the
smooth Lagrangian fibration $\bar f: \bar Z \rightarrow \Gamma$ has
a unique smooth fibre-preserving anti-symplectic involution $\bar
\iota: \bar Z \rightarrow \bar Z$ fixing $\sigma|_{\Gamma}$. Can
$\bar \iota$ be extended to $X$? We have the following result:

\begin{thm}\label{thm:stitched:inv}
 A stitched fibration $\mathcal F = (X,B,f,\sigma, \gamma)$ with invariants $(\bar Z, \ell)$ has a unique smooth fibre-preserving anti-symplectic involution $\iota: X \rightarrow X$ fixing $\sigma$ if and only if
\begin{equation} \label{inv:cond}
 \bar \iota^\ast \ell_k = - \ell_k,
\end{equation}
for every $\ell_k \in \ell$.
\end{thm}
\begin{proof}
Observe that if $\tilde f^+$ and $\tilde f^-$ are smooth, proper,
Lagrangian extensions of $f^+$ and $f^-$ defined on open subsets
$\tilde X^+$ and $\tilde X^-$ of $X$ such that $X^{\pm} \subseteq
\tilde X^{\pm}$, then there are unique fibre-preserving
anti-symplectic involutions $\iota^+: \tilde X^+ \rightarrow \tilde
X^+$ and
 $\iota^-: \tilde X^- \rightarrow \tilde X^-$ fixing $\sigma$. Therefore we may
define $\iota: X \rightarrow X$ to be such that $\iota|_{X^{\pm}} =
\iota^{\pm}$. The question is if $\iota$ is smooth.

Let $\mathcal F_u=(Y, B_u, f_u, \sigma_0, \gamma_0 )$ be a normal
form for $\mathcal F$ (as described in Example~\ref{ex:nf}), and let
$(b,y)$ be the canonical coordinates on $T^*B_u$, then we have the
smooth anti-symplectic involution $\iota_0$ such that $\iota_0(b,y)
= (b, -y)$. We now show that if (\ref{inv:cond}) holds, we can
construct $\mathcal F_u$ so that $\iota_0$ is fibre-preserving with
respect to $f_u$.

We use the same notation as in the construction after
Example~\ref{ex:nf}. Observe that given the sequence $\ell=
\{\ell_k\}_{k\in\N}$ satisfying (\ref{inv:cond}), we can construct a
one parameter family of fiberwise closed sections $\ell(r)$ of
$\mathfrak L^*$ such that (\ref{l:tay}) holds and such that
\begin{equation} \label{ell:cond}
 \bar \iota^\ast \ell(r) = - \ell(r)
\end{equation}
for every small $r$. This can be done by refining the methods used
in \cite{CB-M-stitched} to construct $\ell(r)$ satisfying
(\ref{l:tay}). From $\ell(r)$ we construct $(V,u) \in \mathscr U_Z$
so that the fibres of $u$ are the images $L_{r,b}$ of the graph of
the one form $r dy_1 + \ell(r,b)$ via the symplectomorphism
(\ref{fb:symp}). It can be easily verified that condition
(\ref{ell:cond}) implies that $\iota_0(L_{r,b}) = L_{r,b}$, i.e.
that $\iota_0$ is fibre-preserving with respect to $u$.
 By uniqueness, $\iota_0$ must coincide with the map $\iota$ constructed above under the identification of $\mathcal F_u$ with $\mathcal F$ given by Theorem~\ref{thm: grosso}, part (i). Therefore
$\iota$ is smooth.

Viceversa, suppose now that $\iota$ is smooth, we show that
(\ref{inv:cond}) must hold. Let $\mathcal F_u$ be a normal form for
$\mathcal F$ and let $(V,u) \in \mathscr U_Z$ be the pair defining
$\mathcal F_u$. Then we have that
\begin{equation} \label{i:minus}
  \iota|_{Y^+}: (b,y) \mapsto (b, -y).
\end{equation}
Moreover $u$ satisfies
\begin{equation} \label{i:y+}
  (u \circ \iota)|_{Y^-} = u|_{Y^-}.
\end{equation}
and so also the Taylor expansions with respect to $b_1$ evaluated at
$b_1 = 0$ of the two sides of the above identity must coincide. This
provides a certain relation which must be satisfied by the
coefficients of the Taylor expansion of $u$. Notice that the
coefficients of the Taylor expansion of $u \circ \iota$ only depend
on the coefficients of the Taylor expansions of $u$ and $\iota$, but
the Taylor coefficients of $\iota$ are the same as those of the map
$(b,y) \mapsto (b, -y)$ since (\ref{i:minus}) holds , therefore the
relation among the Taylor coefficients of $u$ implied by equation
(\ref{i:y+}) is the same as the one obtained assuming (\ref{i:y+})
holds with $\iota$ satisfying $\iota(b,y) = (b,-y)$ for all $(b,y)
\in Y$. We now compute this relation in terms of the sequence $\ell
= \{\ell_k\}_{k\in\N}$. Given that $\iota: (b,y) \mapsto (b,-y)$, it
is easy to see $\iota(L_{r,b}) = L_{r,b}$ if and only if the one
parameter family $\ell(r)$ obtained from $u$ satisfies
(\ref{ell:cond}). Therefore the coefficients $\ell_k$ of the Taylor
series (\ref{l:tay}) must satisfy (\ref{inv:cond}).
\end{proof}

\subsubsection{Non proper stitched fibrations}

\medskip
Let $X$ be a smooth symplectic $6$-manifold together with a smooth
Ha\-miltonian $S^1$ action with moment map $\mu:X\rightarrow\R$.
Assume $\mu$ has exactly one critical value $0\in\R$ and a
codimension four submanifold $\Sigma=\Crit \mu$. Let $M$ be a smooth
$2$-dimensional manifold and let $B \subseteq \R \times M$ be a
contractible open neighborhood of a point $(0,m) \in \R \times M$.
Let $\Gamma =B \cap (\{0\}\times M)$. As usual we define $Z =
\mu^{-1}(0)$ and $\bar Z$ the $S^1$ quotient of $Z$ and $X^+=\{
\mu\geq 0\}$, $X^-=\{ \mu\leq 0\}$.

\medskip
We consider fibrations satisfying the following:

\begin{ass}\label{ass: semi-stitched} The map $f:X\rightarrow B$ is a topological $T^3$ fibration with discriminant locus $\Delta \subset \Gamma$ such that $f(\Sigma )=\Delta$ satisfying
\begin{itemize}
\item[(a)] $(X,\omega,f,B)$ is topologically conjugate to a generic singular fibration.
\item[(b)] There is a continuous $S^1$ invariant map $G:X\rightarrow M$ such that
\begin{itemize}
\item[(i)] if $G^{\pm} = G|_{X^{\pm}}$ then $G^+$ and $G^-$ are restrictions of $C^\infty$ maps on $X$;
\item[(ii)] $f$ can be written as $f =  (\mu, G)$ and $f$ restricted to $X^{\pm}$ is a proper map with connected Lagrangian fibres.
\end{itemize}
\item[(c)] There is a connected, $S^1$ invariant, open neighborhood $\mathfrak U\subseteq X$ of $\Sigma$ such that $f(\mathfrak U)=B$ and such that $f_\mathfrak{U}=f|_\mathfrak{U}$ is a $C^\infty$ map with non degenerate singular
points.

\end{itemize}

\end{ass}
This kind of fibrations are studied in \cite{CB-M}. Examples of
fibrations satisfying the above properties can be obtained from the
fibration as in Example \ref{ex amoebous fibr}, after a suitable
perturbation of $f$ near the portion of $\Sigma$ projecting onto the
codimension two part of $\Delta$ (we will recall this smoothing in
the next section). Clearly, the piecewise smoothness occurs along
cylindrical portions of fibres contained in $\mu^{-1}(0)$. Now we
recall some of the basic facts of fibrations satisfying Assumption
\ref{ass: semi-stitched}.

One can construct fibrations of this type as follows. Over $B = D
\times (0,1)$ consider periods given by

\[
\begin{array}{l}
\lambda_{1} = 2\pi db_1, \\
\lambda_2 =dH + \lambda_0, \\
\lambda_3 = db_3,
\end{array}
\]
where $H$ is a smooth function and $\lambda_0=\arg
(b_1+ib_2)db_1+\log |b_1+ib_2|db_2$. If $\Lambda_H$ denotes the
lattice generated by these periods, let  $X^\#=T^\ast B
\slash\Lambda_H$ and denote by $\pi^{\#}$ the projection. Now, in \S
~\ref{generic:singular} we argued that the map $\pi^\#$ can be
extended to a proper map $\pi$ to give a smooth proper Lagrangian
fibration of generic-singular type. This can be achieved by gluing
$\mathfrak U$ to $X^\#$. The moment map of the $S^1$ action is, as
usual, $b_1$. Let $X^{\pm}$, $B^{\pm}$, $Z$ and $\bar Z$ be defined
as usual. Now let $\mathfrak{W}$ and $\mathfrak{U}$ be open $S^1$
invariant neighborhoods of the critical set, such that
$\bar{\mathfrak{W}} \subseteq \mathfrak{U}$. Then, $X^{\circ} =
X-\overline{\mathfrak{W}}$ can be viewed as an open neighborhood of
the zero section of $X^\#$ over which the restriction $\pi^{\circ}$
of $\pi^{\#}$ is a (topologically trivial) Lagrangian open cylinder
fibration (the fibres are homeomorphic to $T^2 \times \R$). The set
$\mathfrak{U} - \bar{\mathfrak{W}}$ covers the two ends of each
fibre. Suppose $u: X^{\circ} \rightarrow B$ is another Lagrangian
open cylinder fibration, whose fibres coincide with the fibres of
$\pi^{\circ}$ over $\mathfrak{U} - \bar{\mathfrak{W}}$, thus the
fibres of $u$ are compactly supported perturbations of the fibres of
$\pi^{\circ}$. If we also assume that $u|_{Z} = \pi^{\circ}$ then we
can define:
\begin{equation}
       f^\circ_{u} = \begin{cases}
               \pi^\circ \quad\text{on} \  X^+, \\
           u \quad\text{on} \ X^-.
       \end{cases}
\end{equation}

The map $f^\circ_{u}$ defines a piece-wise smooth Lagrangian open
cylinder fibration whose fibres coincide with those of $\pi^{\circ}$
on $\mathfrak{U} - \bar{\mathfrak{W}}$. We can therefore glue back
the critical set and define the following proper piecewise smooth
Lagrangian fibration:
\begin{equation} \label{u:st_o}
       f_{u,H}  = \begin{cases}
                \pi \quad\text{on} \  \mathfrak U, \\
           f^\circ_u \quad\text{on} \ X^{\circ}.
       \end{cases}
\end{equation}

Clearly $f_{u,H}: X \rightarrow B$ is well defined and satisfies
Assumption~\ref{ass: semi-stitched}. In \cite{CB-M} it is proved
that any fibration satisfying Assumption~\ref{ass: semi-stitched} is
fiberwise symplectomorphic to $f_{u,H}$ for a certain choice of $u$
and $H$ and therefore $f_{u,H}$ defines a normal form.  Moreover the
invariants that classify such fibrations are given by triples
$(Z^{\#}_{H}, \ell, H_{\Delta})$, where $Z^{\#}_{H}$ is the zero
level set of the $S^1$ moment map (restricted to $X^{\#}$),
$H_\Delta$ is the germ of $H$ along the discriminant $\Delta$, and
$\ell$ a sequence of fiberwise closed sections of $\mathcal L^\ast$,
where $\mathcal L=\ker\bar\pi^\#_\ast$. In this case, each $\ell_k$
is a form with compact support inside cylindrical portions of the
fibres. These invariants classify fibrations as in Assumption
\ref{ass: semi-stitched}. For the details we refer the reader to
\cite{CB-M} \S6.

\begin{thm}\label{thm:stitched_2:inv} A piecewise smooth fibration $f:X\rightarrow B$ satisfying Assumption \ref{ass: semi-stitched} with invariants $(Z^{\#}_{H}, \{\ell_k\}, H_{\Delta})$ and a section $\sigma$ has a unique smooth fibre-preserving anti-symplectic involution $\iota: X\rightarrow X$ fixing $\sigma$ if and only if the invariants satisfy $\bar\iota^\ast \ell_k=-\ell_k$.
\begin{proof}
Given a normal form $(X, f_{u,H})$ for a fibration satisfying
Assumption \ref{ass: semi-stitched}, consider the anti-symplectic
involution $\iota: X \rightarrow X$ constructed in
Section~\ref{generic:singular} preserving the fibres of the smooth
fibration $\pi$. Using the same arguments as in the proof of
Proposition \ref{invo:nodal}, one can prove that $\iota$ also
preserves the fibres of $f_{u,H}$ if the invariants satisfy
$\bar\iota^\ast \ell_k=-\ell_k$. Viceversa given invariants
satisfying this identity on can construct a fibration $f_{u,H}$
having these invariants and whose fibres are preserved by $\iota$.
\end{proof}
\end{thm}

\subsection{Negative fibration.}
Let $f:X\rightarrow B$ be the piecewise smooth fibration  in Example
\ref{ex amoebous_th fibr}. Recall that $X\subset\C^3$ and the
construction of $f$ makes use of a choice of symplectomorphism
$\Phi$ as in (\ref{thin:fi}) giving rise to a fibration whose
discriminant locus $\Delta$ is the amoeba of Figure \ref{fig:
amoeba} after its legs are pinched down to a line. We proved this
fibration is invariant under the standard conjugation on $\C^3$. The
fixed locus consists of 5 connected components, two of which are
sections. The section $\sigma$ fixed by the involution is given by
the choice of any of such sections (cf. Remark
\ref{rem:section-neg}).

In \cite{CB-M}Theorem 7.3, the first two authors propose a method to
make the aforementioned $f$ smoother, obtaining examples of
fibrations of negative type.

\begin{defi} \label{lag:neg}
Let \ $X$ \  be a $6$-dimensional \ symplectic manifold and $B
\subseteq \R^3$ an open subset. A piecewise smooth Lagrangian
fibration  $f: X \rightarrow B$ is called a \textit{Lagrangian
negative fibration} if it satisfies the following properties:
\begin{itemize}
\item[(i)] $f:X\rightarrow B$ is topologically conjugate to the fibration of Example \ref{ex amoebous_th fibr}, i.e. they define the same germ;
\item[(ii)] there exists a submanifold with boundary $D \subset B$, homeomorphic to a closed disc
in $\R^2$, such that $\Delta \cap (B - D)$ consists of three one
dimensional disjoint segments (the legs of $\Delta$) and $f$ is
smooth when restricted to $X - f^{-1}(D)$;
\item[(iii)] On $B -(D\cup \Delta)$, the affine structure induced by the fibration map is simple.
\item[(iv)] $f$ has a section $\sigma$ such that $\sigma (B)$ is a smooth Lagrangian submanifold disjoint from the singular locus $\Sigma\subset X$ of $f$.
\end{itemize}
\end{defi}

We now present an abbreviated description of the smoothing process
that leads to the proof of existence of negative fibrations (details
cf. \cite{CB-M}\S 7) and show that the existence of an
anti-symplectic involution survives this process.

Let $f:X\rightarrow B$ be the fibration in Example \ref{ex
amoebous_th fibr} and a section as in Remark \ref{rem:section-neg}.
Recall that the anti-symplectic involution preserving $f$ is just
conjugation. Let $b_1, b_2, b_3$ standard coordinates in the base
$B\subseteq \R^3$. Then $\Delta$ is contained in the plane $b_1=0$.
Let $\Sigma\subset X$ be the critical surface-- i.e. the locus where
vanishing cycles collapse, a pair of pants projecting onto $\Delta$
under $f$. For positive $M\in\R$, let $\Delta_{h,M} =\Delta\cap\{
b_2 \leq - M \}$. For $M$ large enough, $\Delta_{h,M}$ is
one-dimensional-- i.e. the thin part of the horizontal leg --and let
$\Sigma_{h,M}$ be the portion of $\Sigma$ projecting onto
$\Delta_{h,M}$. For the following analysis, it is convenient to use
$S^1$-invariant coordinates, $t=\mu$, $u_1=z_1z_2$ and $u_2=z_3$.
Then $u_1$ and $u_2$ can be thought of as coordinates on each
reduced space $\mu^{-1}(t)\slash S^1$. On a suitable small
neighborhood $N_{h,M}$ of $\Sigma_{h,M}$ the restriction of $f$ to
$N_{h,M}$ can be explicitly written as:
\[
f=(\mu, G_t)
\]
where
\begin{equation} \label{Gt}
    G_t(u_1, u_2) = \left(
        \log |u_2|, \log \left| \frac{u_1}{\sqrt{|t|+\sqrt{t^2+|u_1|^2}}}-1 \right|
                                                     \right).
\end{equation}

Clearly, $f$ fails to be smooth at $t=0$ since $G_t$ does. In
\cite{CB-M} it is shown that one can perturb $f$ on $N_{h,M}$ by
replacing $G_t$ with a map of type:
\[
\tilde G_t = \left(
        \log |u_2|, \log \left| \frac{u_1}{\rho(|u_1|,t, |u_2|^2)}-1 \right|
                                                     \right).
\]
Here $\rho$ is chosen so that $\tilde G_t$ coincides with $G_t$ away
from $N_{h,M}$ and it is smooth on $N_{h,M}$ (details cf.
\cite{CB-M} Lemma 7.4). It is clear that $\tilde G_t$ is invariant
under the involution on $\C^2$, i.e. under conjugation. The
perturbation $\tilde f=(\mu, \tilde G_t)$ of $f$ is therefore
invariant under the standard involution on $\C^3$.

Similarly, one perturbs $f$ along small neighborhoods of $N_{v,M}$,
and $N_{d,M}$ of $\Sigma_{v,M}$ and $\Sigma_{d,M}$ projecting onto
$B_{v,M}$ and $B_{d,M}$ open neighborhoods of the vertical and
diagonal legs, respectively. This produces an involution-invariant
fibration $\tilde f$.

Now the smoothing needs to be extended to $X_{h,M}:=\tilde
f^{-1}(B_{h,M})$. First observe that the restriction of $\tilde f$
to $X_{h,M}$ is a piecewise smooth fibration satisfying the
hypothesis of Theorem \ref{thm:stitched_2:inv}. Since $X_{h,M}$ has
a fibre-preserving involution $\iota$ fixing a section, if $\{\ell_k\}$ are the invariants of $\tilde f$, then $\bar\iota^\ast\ell_k=-\ell_k$.

Now in \cite{CB-M} Lemma 7.6 it is shown that for some positive $m >
M$ there is a neighborhood $B_{h,m}\subset B_{h,M}\cap\{b_2\leq
-m\}$ and a perturbation of $\tilde f$, making it smooth on
$X_{h,m}:=\tilde f^{-1}(B_{h,m})$. This is achieved by perturbing
the invariants $\ell_k$ of $\tilde f$  in such a way that $\ell_k$
vanish identically on $B_{h,m}\cap\{b_1=0\}$. The perturbed
invariants are (with a slight abuse of notation) of the form
$\nu\ell_k$, where $\nu$ is a bump function on
$B_{h,M}\cap\{b_1=0\}$ vanishing identically on
$B_{h,m}\cap\{b_1=0\}$. Since $\nu$ is a function depending only on
coordinates of the base, it is clear that
$\bar\iota^\ast(\nu\ell_k)=-\nu\bar\ell_k$. Therefore the resulting
fibration after this perturbation is invariant under $\iota$ and the
section resulting from this perturbation is fixed by $\iota$.

One may proceed in an analogous way with the other two legs. This
gives a piecewise smooth fibration which is smooth over large open
neighborhoods $B_{h,m}$, $B_{v,m}$, $B_{d,m}$, of $\Delta_{h,m}$,
$\Delta_{v,m}$, $\Delta_{d,m}$, respectively, and a smooth
fibre-preserving anti-symplectic involution defined on the total
space of the fibration.

Finally, to produce a fibration satisfying the properties of
Definition \ref{lag:neg}, one needs to perturb the fibration away
from a (planar) tubular neighborhood $N$ of $\Delta$. Observe that
the complement of $\Delta$ in the plane $\{b_1=0\}$ consists of
three connected components, $\Gamma_c$, $\Gamma_d$, $\Gamma_e$ which
are the walls of three stitched fibrations $f_c$, $f_d$, $f_e$, each
fibration being the restriction of the fibration obtained in the
previous paragraph. If $\ell^c$, $\ell^d$ and $\ell^e$ are the
corresponding invariants, then Theorem \ref{thm:stitched:inv}
implies that $\bar\iota^\ast\ell^c_k=-\ell^c_k$,
$\bar\iota^\ast\ell^d_k=-\ell^d_k$ and
$\bar\iota^\ast\ell^e_k=-\ell^e_k$. Now, in \cite{CB-M} Lemma 7.12,
it is shown that $f_c$, can be made smooth away from
$N\cap\Gamma_c$. As before, this is accomplished after deforming
$\ell^c_k$ to $\tilde\ell^c_k=\rho\ell^c_k$ for a suitably chosen
bump function $\rho$ on $\Gamma_c$. Again, being $\rho$ dependent on
coordinates on the base, implies that the resulting fibration is
still $\iota$-invariant and the resulting section fixed by $\iota$. One proceeds in a similar way with $f_d$
and $f_e$. This completes the required smoothing of $f$.

Observe that if $D$ is the region over which $f$ fails to be smooth,
there are regions $D'\subset D$ and $B'\subset B$ such that
$B'\cap\{b_1=0\}\subset D'$ where the section $\sigma$ obtained
after the smoothing of $f$ remains unchanged, i.e. $\sigma (B')$
coincides with the section in Remark \ref{rem:section-neg}. It also
follows that $\sigma (B)$ is smooth. This completes the proof the
following:

\begin{thm}\label{thm:neg-fib}
Let $f:X\rightarrow B$ be a Lagrangian negative fibration. Then
there is a Lagrangian section $\sigma$ not intersecting the singular
locus $\Sigma\subset X$ of $f$ and unique smooth fibre preserving
anti-symplectic involution $\iota_{f,\sigma}$ of $X$ preserving the
fibres of $f$ and fixing $\sigma$.
\end{thm}

\begin{rem}\label{rem:no-conj}
Let $\sigma_1$ and $\sigma_2$ Lagrangian sections of a negative fibration $f:X\rightarrow B$ and $B'\subset B$ and $D'\subset D$ as above. Since $f$ is stitched along $f^{-1}(D)\subset X$, it follows from Remark \ref{rem:no-trans-stitched}, that, in general, there is no symplectomorphism $t$ of $X$ such that $f\circ t$ and $t\circ\sigma_1=\sigma_2$. This contrasts with the nodal, generic-singular and positive models for which $t$ always exists.
\end{rem}

\begin{lem}\label{lem:conj-neg} Let $\sigma_1$ and $\sigma_2$ be sections of a negative fibration $f:X\rightarrow B$ and $D\subset B\cap\{b_1=0\}$ the locus over which $f$ fails to be smooth. If there exists an open neighborhood $B' \subset B$ of $D$ such that $\sigma_1|_{B'} =\sigma_2|_{B'}$, then there is a unique symplectomorphism $t$ of $X$ such that $f\circ t = f$ and $t\circ\sigma_1=\sigma_2$
\end{lem}
\begin{proof}
On $X_\circ=f^{-1}(B-B')$, the fibration is smooth. Corollary \ref{cor:tw} and Lemma \ref{lem:tw:non-proper-gen} give a unique symplectomorphism $t_\circ$ of $X_\circ$ sending $\sigma_1|_{B-B'}$ to $\sigma_2|_{B-B'}$. Extending $t_\circ$ to $X$ as the identity map on $X-X_\circ$ gives a smooth symplectomorphism $t$ with the required properties.
\end{proof}

\section{Global existence}\label{sec:global}

Let $(B,\Delta, \mathscr A)$ be a compact simple integral affine
manifold with singularities. Let $\mathcal N$ be the set of negative vertices of $\Delta$ and let $(\Delta_{\blacklozenge}, \{
D_{p} \}_{p \in \mathcal N})$ be a localized thickening and let
$(B_\blacklozenge, \Delta_\blacklozenge, \mathscr A_\blacklozenge)$
be the integral affine manifold as in \S \ref{sec:C}. Then there is
a smooth symplectic manifold
\[
X_\blacklozenge =T^\ast B_\blacklozenge\slash \Lambda_\blacklozenge
\]
where $\Lambda_\blacklozenge$ is the period lattice induced by
$\mathscr A|_{B_\blacklozenge}$, and a Lagrangian submersion:
\[
f_\blacklozenge:X_\blacklozenge\rightarrow B_\blacklozenge .
\]
Notice that if $B_0=B-\Delta$, $\Lambda_0$ is the lattice induced by
$\mathscr A$, $X_0=T^\ast B_0\slash\Lambda_0$, and
$f_0:X_0\rightarrow B_0$ the standard projection, then
$X_\blacklozenge\subset X_0$ and
$f_\blacklozenge=f_0|_{X_\blacklozenge}$.

Let $\sigma_0$ be a section of $f_0$ which can be taken to be
induced by the zero section on $T^\ast B_0$. Then, Corollary
\ref{cor:m} implies there is a unique fibre preserving
anti-symplectic involution $\phi_0$ of $X_0$ also preserving
$\sigma_0$. With abuse of notation, denote by $\sigma_0$ and
$\phi_0$ their restrictions to $X_\blacklozenge$ and $B_\blacklozenge$ respectively.
Theorem \ref{class_C} gives a class $\mathcal C$ of fibrations $f:X\rightarrow B$ where  $X$ is the compact symplectic manifold
obtained from $X_\blacklozenge$ after gluing models of generic,
positive and negative fibrations as in \S\ref{sec:local} over
$\Delta$ and matching local sections of each local model with $\sigma_0$. The latter provides the fibration with a section $\sigma$.

\subsection{The class $\Cc$} We now impose extra conditions on the sections of fibrations of class $\CC$.

\begin{defi}\label{defi:class_Cc}
Let $f:X\rightarrow B$ be a fibration of class $\mathcal C$ with a
section $\sigma$ such that $\sigma (B)\cap \Crit f=\varnothing$
where $\Crit f\subset X$ is the singular set of $f$.  Assume that
identifications of neighborhoods of singular fibres with the local
models of Section~\ref{sec:local} are fixed. For each negative
vertex $p\in\mathcal N$, let $B_{p}\subset B$ be a small open
neighborhood of $p$ such that $D_{p}\subset B_{p}$, where $D_p$ is
the locus over which $f$ is piecewise smooth.  Let
$f^-:X^-\rightarrow B^-$ be the model for the negative fibration and
let $\sigma^-$ be a choice of section of $f^-$ fixed by the local
anti-symplectic involution as in Theorem \ref{thm:neg-fib}. We say
that $\sigma$ is of class $\Cc$ if for each $p\in\mathcal N$, the
restriction of $\sigma$ to $B_{p}$ coincides with $\sigma^-$.
\end{defi}

Notice that the definition of $\Cc$ clearly depends on the choice of
$\sigma^-$ and $\{B_p\}_{p\in\mathcal N}$.  Notice also that another
section $\sigma'$ is of class $\Cc$ if and only if $\sigma'$
coincides with $\sigma$ when restricted to each $B_p$.

\subsection{Proof of Theorem \ref{tm:tw}} Let $f:X\rightarrow B$ a fibration of class $\CC$ with two sections $\sigma_1, \sigma_2 \in \Cc$.  Let $X_\blacklozenge\subset X$ and $f_\blacklozenge :X_\blacklozenge\rightarrow B_\blacklozenge$ be the Lagrangian submersion as above. By Corollary \ref{cor:tw} there is a unique fibre-preserving symplectomorphism $t_0: X_\blacklozenge \rightarrow X_\blacklozenge $ sending $\sigma_{1}$ to $\sigma_{2}$. We will show that $t_0$ extends to $X$.

Since $\sigma_1,\sigma_2\in\Cc$, it follows from Definition \ref{defi:class_Cc} that for each negative vertex $p\in\mathcal N$,  $\sigma_1|_{B_p}=\sigma_2|_{B_p}$. Trivially, there is a unique local fibre-preserving symplectomorphism sending the restriction $\sigma_1|_{B_p}$ to the restriction $\sigma_2|_{B_p}$. Similarly, Lemma \ref{lem:tw:non-proper-pos} guarantees that for each positive vertex $v$ of $\Delta$, there is an open neighborhood $B_v\subset B$  of $v$ and a unique local fibre-preserving symplectomorphism sending the restriction $\sigma_1|_{B_v}$ to the restriction $\sigma_2|_{B_v}$. For the edges of $\Delta$ one applies Lemma \ref{lem:tw:non-proper-gen} analogously. Each of these local symplectomorphisms provide a local extension of $t_0$ to $X$. By uniqueness, these extensions glue together along common intersection, giving a unique extension $t$ of $t_0$ to $X$. Details are left to the reader.

\begin{rem}
Notice that due to the piecewise smoothness of $f:X\rightarrow B$ of class $\mathcal C$, if $\sigma_1$ is of class $\Cc$ but $\sigma_2$ is not, it cannot follow that there is a symplectomorphism $t$ of $X$ such that $f\circ t=f$ and $t\circ\sigma_1=\sigma_2$ (cf. Remark \ref{rem:no-conj}).
\end{rem}

\subsection{Proof of Theorem \ref{tm:m}}
It is enough to find a fiber-preserving anti-symplectic involution
$\phi$ fixing one section $\sigma'\in\Cc$. In fact, if $\sigma$ is
any other section in $\Cc$ and $t$ is the symplectomorphism taking
$\sigma$ to $\sigma'$ constructed in Theorem~\ref{tm:tw}, then
$\phi_{f,\sigma} = t^{-1} \circ \phi \circ t$ is the anti-symplectic
involution fixing $\sigma$.

Consider as above, the anti-symplectic involution $\phi_0$ of
$X_{\blacklozenge}$ fixing the section $\sigma_0$. We need to show
that the section $\sigma_0$ extends to a section $\sigma'\in\Cc$ and
the involution $\phi_0$ of $X_0$ extends to a smooth
fibre-preserving anti-symplectic involution $ \phi $ of $X$ fixing
$\sigma'$. The proof follows immediately from Theorem \ref{class_C}
and the results of \S\ref{sec:local}. Let us denote by
$f^\nu:X^\nu\rightarrow B^\nu$ a fibration of either
generic-singular, positive or negative type, used in the
compactification as in Theorem \ref{class_C}.  This presumes that
each affine base, $B_0^\nu=B^\nu-\Delta^\nu$, is locally affine
isomorphic to $U_0=U-\Delta\cap U$, where $U$ is a suitable
neighborhood of $x\in\Delta$, and $x$ is either an edge point, a
positive or a negative vertex. If $\nu$ is either generic or
positive, we let $\sigma^\nu$ be any choice of a section of $f^\nu$
fixed by the local anti-symplectic involution, not intersecting the
critical locus of $f^\nu$. If $\nu$ is a negative vertex, the choice
of $\sigma^\nu$ is $\sigma^-$ as in Definition \ref{defi:class_Cc}.

 Then, the affine isomorphism induces a symplectomorphism of
bundles $\Phi^\nu$ and a commuting diagram
\begin{equation}\label{eq:glue}
\begin{CD}
X^\nu_0  @>\Phi^\nu>> f^{-1}_\blacklozenge (U_0)\\
@Vf^\nu VV  @VVf_\blacklozenge V\\
B^\nu_0 @>A^\nu>> U_0
\end{CD}
\end{equation}
where $\Phi^\nu (\sigma^\nu|_{B^\nu_0})=\sigma_0\circ A^\nu$ and
$A^\nu$ extends continuously to $B^\nu$. Then $X^\nu$ is glued to
$X_\blacklozenge$ over $U$ using $\Phi^\nu$. Moreover, this gluing
extends $\sigma_0$ to a smooth section on $U$. The gluing of two
generic-singular fibrations along common edges in $\Delta$ requires
taking care of further technicalities, as it involves gluing along
singular fibres. A smooth symplectic deformation of the fibrations
along a common intersection may be required but, in any case, two
generic-singular fibrations can be glued matching its corresponding
prescribed Lagrangian sections (cf. \cite{CB-M} Proposition 4.18).

Now $f_\blacklozenge$ and each $f_\nu$ carry a unique fibre
preserving smooth anti-symplectic involution $\phi_0$ and $\phi^\nu$
fixing $\sigma_0$ and $\sigma^\nu$, respectively. Since $\sigma_0$
and $\sigma^\nu$ coincide over $U_0$, it follows that $\phi_0$ and
$\phi^\nu$ coincide along $f^{-1}_\blacklozenge (U_0)$. Then
$\phi_0$ extends smoothly to $f^{-1}_\blacklozenge(U)$. Repeating
this process for a suitable open cover $\{U\}$ of $\Delta$ produces
the required section $\sigma'$ of $f$, and the extension $\phi$ of $\phi_0$. This completes the proof of Theorem \ref{tm:m}. By construction $\sigma\in\Cc$.

\section{Examples}
The same arguments discussed in \S\ref{sec:global} apply in
dimension $n=2$:

\begin{thm}\label{thm:2dim} Let $(B,\Delta,\mathscr A)$ be a 2-dimensional simple affine manifold with singularities and let $f:X_0\rightarrow B_0=B-\Delta$ be the Lagrangian submersion of $X_0=T^\ast B_0\slash\Lambda$ onto $B_0$. Then
\begin{itemize}
\item[(i)] There is a symplectic manifold $X$ and a Lagrangian fibration $f:X\rightarrow B$ such that $f|_{X_0}=f_0$;
\item [(ii)]If a Lagrangian section $\sigma$ is specified which avoids the critical points, then there is a unique fibre preserving anti-symplectic involution $\phi_{f,\sigma}$ fixing $\sigma$.
\item[(iii)]If two Lagrangian sections $\sigma_1$ and $\sigma_2$ are specified (both avoiding the critical points), there is a unique symplectomorphism $t:X\rightarrow X$ such that $f\circ t=f$ and $t\circ \sigma_1=\sigma_2$.
\end{itemize}
\end{thm}
The first claim is the content of \cite{CB-M} Theorem 3.22, while
the second and third claims are new. The proof of (ii) is a verbatim
of the one in dimension 3, where one can use the model for a
focus-focus fibration of Section~\ref{smooth:ff} together with the
given anti-symplectic involution. The proof of (iii) is the same as
the proof of Theorem \ref{tm:tw}.

\subsection{The K3}
Starting with explicit examples of integral affine base $B\cong
S^2$ with 24 singularities Leung and Symington \cite{LeungSym} illustrate how part (i) of Theorem \ref{thm:2dim} can be used to build well known Lagrangian fibrations on a
symplectic 4-manifold $X\cong K3$ with a section (see also \cite{CB-M}, Example 3.16).
The construction involves making several choices, which produce different germs of
Lagrangian fibrations. So even though the compactification $X$ is
the same (modulo symplectomorphism) regardless of the choices made,
there are actually infinitely many germs of Lagrangian fibrations
with the same topology (cf. \cite{CB-M} Corollary 3.24). Given a
choice of such fibration germ $f$, part (ii) gives a unique fibre
preserving anti-symplectic involution.

In this case, the fixed locus of $\phi_f$ is a Lagrangian
submanifold with 2 connected components: one of them is a sphere
(i.e. the section) and the other is a genus $g=10$ surface $\Sigma$
which is a 2:1 branch cover of $S^2$, with 24 branch points.

\subsection{Almost toric 4-manifolds}
Symington and Leung \cite{LeungSym} propose a class of symplectic
4-manifolds with Lagrangian fibrations having focus-focus and toric
singular fibres, called \textit{almost toric}. Within this class,
the integral affine bases which are simple (\textit{simple} in the
sense of Theorem \ref{thm:2dim}) are the disc $D^2$, the cylinder
$S^1\times I$, the Klein bottle, the sphere $S^2$ and $\R\PP^2$.
Namely, these are the only cases that have singularities of nodal
type. Theorem \ref{thm:2dim} equips each of the corresponding
fibrations with fibre-preserving anti-symplectic involutions. For
instance, the Enriques surface is equipped with a Lagrangian
fibration over $\R\PP^2$ with 12 focus-focus singularities and a
fibre preserving anti-symplectic involution. The base $S^2$ gives a
K3 surface discussed above.

\subsection{The quintic}

Starting with an explicit example of affine 3-manifold with
singularities proposed by Gross \cite{Gross_SYZ_rev} Example 4.3,
the first two authors use Theorem \ref{class_C} to produce a
symplectic 6-manifold $X$ homeomorphic to a smooth quintic 3-fold.
Now Theorem \ref{tm:m} shows that such manifold has a fibre
preserving anti-symplectic involution.

\subsection{Mirror pairs}

The example above generalizes to a much wider class. When $B$ is an
integral affine 3-manifold arising from toric degeneration in the
sense of Gross and Siebert, Theorem \ref{class_C} produces pairs of
SYZ dual Lagrangian fibrations, with total spaces homeomorphic to
mirror pairs of Calabi-Yau manifolds (cf. \cite{CB-M} for details).
Now Theorem \ref{tm:m} equips these pairs of symplectic 6-manifolds
with fibre-preserving anti-symplectic involutions fixing a section.

\medskip
In the examples discussed above, the fixed locus set $\Sigma$
appears to have nice topological properties. For instance, for $X$
Calabi-Yau, there is an intriguing relation between the mod 2
cohomology of $\Sigma$ and the Hodge numbers of $X$. These
properties are being further investigated in \cite{CB-M-top}.

\section{Fiber-preserving anti-symplectomorphisms}\label{sec:as}
In this section, we prove Proposition~\ref{pr:ac}.

In the following lemmas, $f: X \rightarrow B$ is a Lagrangian
fibration that is a smooth submersion, and $\sigma$ is a smooth
Lagrangian section of $f.$ We denote by $\Lambda$ the lattice bundle
and by $\Theta : T^*B/\Lambda \rightarrow X$ the symplectomorphism
of Proposition~\ref{pr:aa} applied to $\sigma.$ We denote by $Z$ the
zero sections of $T^*B$ and $T^*B/\Lambda,$ and we denote by $\pi$
the canonical projections to $B.$ We denote by $-\id$ the
anti-symplectomorphisms of $T^*B$ and $T^*B/\Lambda$ given by
negative the identity map on each fiber.

Let $\eta$ be a $1$-form on $B.$ We define a symplectomorphism
$T_\eta : T^*B \rightarrow T^*B$ by
\[
T_\eta(p,\xi) = (p, \xi + \eta(p)), \qquad \forall p \in B,\; \xi
\in T^*_pB.
\]
We also denote by $T_\eta$ the symplectomorphism that $T_\eta$
induces on $T^*B/\Lambda.$

\begin{lem}\label{lem:phi}
Assume that $\pi_1(B) = \{1\}.$ Let $\phi$ be an
anti-symplectomorphism of $X$ such that $f \circ \phi = \phi.$ Then
$\phi^2 = \id_X.$ In particular,
\begin{equation}\label{eq:-id}
\phi = \Theta \circ T_\eta \circ (-\id) \circ \Theta^{-1}.
\end{equation}
\end{lem}
\begin{proof}
Define
\[
Z' = \Theta^{-1} \circ \phi \circ \Theta \circ Z.
\]
It is easy to see that $Z'$ is a Lagrangian section of
$T^*B/\Lambda.$ Since $\pi_1(B) = \{1\},$ we may lift $Z'$ to a
Lagrangian section $\widetilde Z'$ of $T^*B.$ Let $\eta$ be the one
form on $B$ such that $\widetilde Z'$ is its graph. Clearly, $\pi
\circ T_\eta = \pi$ and $T_\eta \circ Z = Z'.$ By the uniqueness
claim of Corollary~\ref{cor:tw} applied to the Lagrangian fibration
$\pi : T^*B/\Lambda \rightarrow B,$ we conclude that
\[
T_\eta = \Theta^{-1} \circ \phi \circ \Theta \circ (-\id).
\]
Formula \eqref{eq:-id} follows. Observe that
\[
T_\eta \circ (-\id) = (-\id) \circ T_{-\eta}, \qquad T_{-\eta} =
T_\eta^{-1}.
\]
Consequently, $(T_\eta \circ (-\id))^2 = \id_{T^*B/\Lambda}.$ The
lemma follows.
\end{proof}

We omit the proof of the following lemma since it is similar and we
do not use it.

\begin{lem}\label{lem:t}
Assume that $\pi_1(B) = \{1\}.$ Let $t$ be a symplectomorphism of
$X$ such that $f \circ t = f.$ There exists a $1$-form on $B$ such
that
\[
t = \Theta \circ T_\eta \circ \Theta^{-1}.
\]
\end{lem}

\begin{proof}[Proof of Proposition~\ref{pr:ac}]
Since smooth fibers are dense and the claim is a closed condition,
we may assume without loss of generality that $f: X\rightarrow B$ is
a smooth submersion. Since $\phi_f$ preserves fibers of $f,$ the
claim is local on the base $B.$ So, without loss of generality we
focus on the special case when $B$ is the $n$-disk. Equation
\eqref{eq:sq} follows from Lemma~\ref{lem:phi}. Equation
\eqref{eq:ac} follows formally from equation \eqref{eq:sq}. Indeed,
$\phi_f \circ t$ is an anti-symplectomorphism such that $f\circ
\phi_f \circ t = f.$ So, we conclude $\phi_f \circ t \circ \phi_f
\circ t = \id_X,$ which implies equation \eqref{eq:ac}.
\end{proof}

\section{Gradings}\label{sec:g}
In this section we will explain the definition of the grading of a
Lagrangian submanifold $L \subset X$ in the special case where $X$
is a symplectic Calabi-Yau manifold. Then we will assume that $f: X
\rightarrow B$ is a special Lagrangian fibration and
$\phi_{f,\sigma}$ is anti-holomorphic as well as anti-symplectic. In
this case, we conclude that $\I_{f,\sigma}$ shifts the natural
grading on the fibers of $f$ by $\dim_\C X.$

The notion of a grading for a Lagrangian submanifold was introduced
by Kontsevich \cite{Kontsevich}. Here we follow a slightly modified
version of the exposition of \cite{TY}. We use the generalized
definition of special Lagrangian submanifolds due to Salur \cite{Sa}
that applies to symplectic Calabi-Yau manifolds that may not have an
integrable complex structure.

Let $X$ be a symplectic manifold with symplectic form $\omega.$ An
almost complex structure on $J$ on $X$ is said to be $\omega$-tame
if
\[
\omega(\xi,J\xi) > 0
\]
for all $\xi \neq 0.$ We define the first Chern class $c_1(TX)$ by
choosing an $\omega$-tame almost complex structure on $X.$ The
definition of $c_1$ only depends on $\omega$ because the space of
$\omega$-tame almost complex structures is contractible.

>From now on we assume that $(X,\omega)$ is a symplectic
$2n$-real-dimensional Calabi-Yau manifold, i.e that $c_1(TX) = 0.$
We fix an $\omega$-tame almost complex structure $J$ on $X.$ Use $J$
to decompose complex valued differential forms on $X$ by type. Fix a
nowhere vanishing $(n,0)$-form $\Omega$ on $X.$ The existence of
$\Omega$ is guaranteed by the Calabi-Yau condition. We emphasize
that we do not require $\Omega$ to be closed. Finally, we define the
metric $g$ by
\[
g(\xi,\eta) = \frac{\omega(\xi,J\eta) + \omega(\eta,J\xi)}{2}.
\]

Let $L \subset X$ be a Lagrangian submanifold. A small
generalization of arguments of \cite{HL} shows that
\begin{equation}\label{eq:oml}
\Omega|_L = \psi e^{\pi i\theta} vol_g,
\end{equation}
where $\psi$ is a strictly positive real-valued function, $\theta$
is an $S^1$-valued function, and $vol_g$ is the volume form of $L$
induced by $g.$ If the Maslov class of $L$ vanishes, then $\theta$
can be lifted to a real valued function. The choice of a real-valued
lift of $\theta,$ which we also denote by $\theta,$ is a grading of
$L.$ A graded Lagrangian submanifold is called special Lagrangian if
the grading $\theta$ is constant.

For the rest of this section, we assume that $f: X \rightarrow B$ is
a special Lagrangian fibration. That is, each fiber of $f$ contains
a relatively open dense subset that is a smooth special Lagrangian
submanifold of $X.$ We assume that $f$ has a section $\sigma,$ and
we assume that $X$ has an anti-symplectic involution $\phi$
satisfying conditions \eqref{eq:2c}. We assume also that $\phi$ is
anti-$J$-holomorphic and
\begin{equation}\label{eq:om}
\phi^* \Omega = \overline\Omega.
\end{equation}
In Lemma~\ref{lm:om} below, we show that assumption \eqref{eq:om} is
not hard to satisfy given the previous assumptions. Moreover, we
have the following lemma.
\begin{lem}\cite{CMS2}
Let the Lagrangian fibration $f: X\rightarrow B$ be a smooth
submersion. For each $J,$ there exists a choice of $\Omega$ such
that $f$ is special Lagrangian.
\end{lem}

Let $\theta_\sigma$ denote a grading on the Lagrangian submanifold
given by the section $\sigma$ and let $\theta_y$ denote a grading on
the fiber $L_y$ of $f.$
\begin{lem}
The gradings $\theta_\sigma$ and $\theta_y$ must satisfy
\[
\theta_\sigma \in \Z, \qquad \theta_y \in n/2 + \Z.
\]
\end{lem}
\begin{proof}
According to Corollary~\ref{cor:m}, $\phi$ acts on each smooth fiber
$L_y$ of $f$ by a diffeomorphism of sign $(-1)^n.$ So, equating the
phase on each side of \eqref{eq:om} and using the fact that
$\theta_y$ is constant on $L_y,$ we have
\[
e^{\pi i\theta_y} = (-1)^n e^{-\pi i\theta_y}
\]
We conclude that $\theta_y  \in n/2 + \Z.$ On the other hand, $\phi$
acts on $\sigma$ by the identity map. So, the same argument implies
that $\theta_\sigma \in \Z.$
\end{proof}

As noted previously, the mirror correspondence maps $\sigma$ along
with the appropriate local system, spin structure and grading to the
structure sheaf $\O_Y.$ We would like to identify the choice of
grading $\theta_\sigma$ that corresponds to $\O_Y.$ Since $\O_Y$ is
fixed under $\D,$ for consistency of Conjecture~\ref{cj:}, we must
assume that $\I_{f,\sigma}(\sigma,\theta_\sigma) =
(\sigma,\theta_\sigma).$ It follows that
\begin{equation*}
\theta_\sigma = -\theta_\sigma = 0.
\end{equation*}

To fully determine the choice of $\theta_y$ that makes $L_y$ into
the mirror of $\sky_y,$ we employ the mirror correspondence once
again. Since $\sigma$ is sent by the mirror correspondence to the
structure sheaf $\O_Y,$ we should have an isomorphism of graded
vector spaces
\[
m_{f,\sigma} : HF^*(\sigma,L_y) \ris \rsh(\O_Y,\sky_y) \simeq \C,
\]
where the grading of $\C$ is $0.$ We will deduce $\theta_y$ from the
definition of the grading on $HF^*(\sigma,L_y).$

We recall the definition of the grading on $HF^*.$ Let
$L_1,L_2\subset X,$ be two transversely intersecting graded
Lagrangian submanifolds with gradings $\theta_1,\theta_2.$ By
definition, $HF^*(L_1,L_2)$ is the cohomology of the complex
$CF^*(L_1,L_2),$ which is generated by the intersection points of
$L_1$ and $L_2.$ The grading of a point $p\in L_1 \cap L_2$ is
defined as follows. Identify $T_pX$ with $\C^n$ by a complex linear
transformation $t$ taking $L_1$ to $\R^n \subset \C^n$ and $L_2$ to
$M \cdot \R^n.$ Take $M$ to be unitary. So, it is conjugate to a
diagonal matrix of the form
\begin{equation*}
M = \begin{pmatrix}e^{i\pi\alpha_1} & 0 & 0 & \cdots  &  0 \\0 &
e^{i\pi\alpha_2} & 0 & \cdots &   0 \\ 0 & 0 & e^{i\pi\alpha_3} & &
\vdots \\\vdots & \vdots & \vdots      & \ddots &   0 \\ 0 & 0 & 0 &
0 & e^{i\pi\alpha_n} \\\end{pmatrix}
\end{equation*}
where $\alpha_i \in (0,1).$ Set $\alpha = \sum_i \alpha_i.$ Define
the grading of $p$ to be
\begin{equation}\label{eq:gr}
ind_p(L_1,L_2) = \alpha - \theta_2(p) + \theta_1(p).
\end{equation}

\begin{lem}
Let $L_y$ be a smooth fiber. Assuming $\theta_\sigma \equiv 0,$ and
$HF^*(\sigma,L_y)$ is a one dimensional vector space of grading $0,$
it follows that $\theta_y = n/2.$
\end{lem}

\begin{proof}
By definition of a section, there is a unique intersection point $p
\in \sigma \cap L_y.$ Let $t: T_p X \rightarrow \C^n$ such that
$t(T_p\sigma) = \R^n.$ It follows from Corollary~\ref{cor:m} that we
can choose $t$ so that $t(L_y) = i \R^n.$ Then for $i = 1,\cdots,n,$
we have $\alpha_i = 1/2.$ So, $\alpha = n/2.$ Rearranging equation
\eqref{eq:gr} we obtain
\[
\theta_y(p) = \alpha - ind_p(\sigma,L_y) + \theta_\sigma(p)  =
\alpha = n/2.
\]
\end{proof}
It follows that $\I_{f,\sigma}$ shifts $\theta_y$ by $n = \dim_\C
X.$

\begin{rem}
The fact that the natural grading for a torus fiber is $n/2$ has
been observed previously by Douglas in the context of
$\Pi$-stability \cite{Do}. See also \cite{Aroux-SYZ}.
\end{rem}

We close this section with a lemma that shows that assumption
\eqref{eq:om} follows from the other assumptions under mild
conditions.
\begin{lem}\label{lm:om}
Let $(X,\omega)$ be a symplectic Calabi-Yau manifold with
$\omega$-tame almost complex structure $J$ and nowhere-vanishing
$(n,0)$-form $\Omega.$ Let $\phi$ be an anti-symplectic involution
of $X$ that is also anti-$J$-holomorphic. If $X$ is simply connected
or if $\Omega$ is closed, then there exists a smooth complex valued
function $g$ on $X$ such that $\widehat\Omega = g\Omega$ satisfies
condition \eqref{eq:om}. If $f : X \rightarrow B$ is a special
Lagrangian fibration with respect to $\Omega$ and $f \circ \phi =
f,$ then $f$ is also a special Lagrangian fibration with respect to
$\widehat \Omega.$
\end{lem}
\begin{proof}
Since $\Lambda^{0,3}(T^*X)$ is a line bundle, there exists a complex
valued function $h$ such that $\phi^* \Omega = h \overline{\Omega}.$
It follows from the fact that $\phi$ is an involution that
\[
h\circ\phi = h^{-1}.
\]
So, if $h$ has a square root $h^{1/2}$ we can take $g = h^{1/2}.$
Clearly, if $X$ is simply connected, then $g$ has a square-root.
Alternatively, if $\Omega$ is closed, then $J$ is integrable and
$\Omega$ is holomorphic \cite{Hitchin}. So, both $\phi^*\Omega$ and
$\overline \Omega$ are anti-holomorphic and therefore so is $h.$ It
follows that $h$ is constant and hence has a square-root.

To prove the final claim, we show that the phase of $h$ is constant
on fibers of $f.$ It follows that if $f$ is special Lagrangian with
respect to $\Omega$ then it is also special Lagrangian with respect
to $\widehat \Omega = h^{1/2}\Omega.$ Indeed, let $L_y$ be a fiber
of $f.$ In the notation of equation \eqref{eq:oml}, using the fact
that $\theta_y$ is constant, we have
\[
(-1)^n(\psi_y\circ \phi) e^{i\pi\theta_y} vol_g = \phi^*
\Omega|_{L_y} = h \overline\Omega|_{L_y} = h|_{L_y} \psi_y
e^{-i\pi\theta_y} vol_g.
\]
It follows that
\[
h|_{L_y} = (-1)^n(\psi_y \circ \phi) \psi_y^{-1} e^{2\pi i
\theta_y},
\]
which has constant phase.
\end{proof}

\section{Coherent Sheaves}\label{sec:cs}
In this section, we prove Theorems~\ref{tm:o1} and~\ref{tm:o2}. The
proof uses Theorem~\ref{thm:or} below, which was proven by D. Orlov
\cite{Or}.

Let $M$ and $X$ be smooth projective varieties over a field $k.$ For
any object $E$ of $\DC(M \times X)$ we can define a functor
\[
\Phi_E : \DC(M) \rightarrow \DC(X)
\]
as follows. Let $p : M \times X \rightarrow M$ and $\pi : M\times X
\rightarrow X$ denote the projections to $M$ and $X$ respectively.
Define
\begin{equation}\label{eq:E}
\Phi_E (\bullet) = R\pi_*(E \overset{L}{\otimes} p^*(\bullet)).
\end{equation}
\begin{thm}\label{thm:or}
Let $F$ be an exact functor from $\DC(M)$ to $\DC(X),$ where $M$ and
$X$ are smooth projective varieties. Suppose $F$ is full and
faithful and has a right (and, consequently, a left) adjoint
functor. Then there exists an object $E$ of $\DC(M \times X)$ such
that $F$ is isomorphic to the functor $\Phi_E$ defined by the
rule~\eqref{eq:E}, and this object is unique up to isomorphism.
\end{thm}

The following corollary parallels Proposition~\ref{pr:aa}. As above,
$\sky_y$ denotes the skyscraper sheaf at a point $y \in Y.$

\begin{cor}~\label{cor:maa}
Let $Y$ be a smooth projective variety. Let $F : \DC(Y) \rightarrow
\DC(Y)$ be an auto-equivalence such that
\begin{equation*}
F(\sky_y) \simeq \sky_y, \;\;\forall y \in Y, \qquad \qquad F(\O_Y)
\simeq \O_Y.
\end{equation*}
Then, $F$ is isomorphic to the identity functor.
\end{cor}
\begin{proof}[Proof of Corollary~\ref{cor:maa}]
We apply Theorem~\ref{thm:or} in the case that $M = X = Y.$ Let $E$
be the object of $\DC(Y \times Y)$ associated to $F$ by
Theorem~\ref{thm:or}. Let $\Delta : Y \rightarrow Y \times Y$ denote
the diagonal map. We also use $\Delta$ to denote the diagonal
subvariety.

Given an object $C$ of $\DC(Y)$ we define $\supp(C)$ to be the
closed subset of the underlying topological space of $Y$ that is the
union of the supports of all the cohomology sheaves of $C.$

First, we prove that $\supp(E) \subset \Delta.$ Indeed, define
\[
i_y : Y \rightarrow Y \times Y
\]
to be the inclusion of $y \times Y$ into $Y \times Y.$ Since
\[
p^* \sky_y \simeq \O_{y \times Y}
\]
we have
\[
E \overset{L}{\otimes} p^*\sky_y \simeq E|_{y \times Y},
\]
where the restriction is in the derived sense. Since $\pi|_{y \times
Y}$ is the identity map,
\[
F(\sky_y) = R\pi_*(E \overset{L}{\otimes} p^*\sky_y) \simeq Li_y^*
E.
\]
So, by assumption,
\begin{equation}\label{eq:fib}
Li_y^* E \simeq \sky_y.
\end{equation}
So, the only fibers of $E$ which are not zero are on $\Delta.$

Next, we prove that $E$ can be represented by a complex concentrated
in degree zero. Since $\supp(E) \subset \Delta,$ we know that $E$ is
the push-forward of a complex of sheaves supported on the diagonal
with some possibly non-reduced scheme structure. So, it suffices to
work in a neighborhood of the diagonal. Locally, in a neighborhood
of the diagonal, $\pi$ is affine. Indeed, if $U \subset Y$ is an
open affine, then
\[
\Delta \cap \pi^{-1}(U) \subset U \times U \subset \pi^{-1}(U).
\]
We call this property (LA). Since
\[
p^* \O_Y = \O_{Y \times Y},
\]
we have by assumption,
\[
R\pi_*(E) = F(\O_Y) \simeq \O_Y.
\]
But, by property (LA), $R\pi_*$ coincides with $\pi_*.$ Moreover,
affine push-forward cannot send a non-trivial sheaf to zero.
Therefore, like its push-forward, $E$ must be concentrated in degree
$0.$

Next, we construct an isomorphism
\begin{equation}\label{eq:m1}
\O_\Delta \ris E|_{\Delta}.
\end{equation}
Indeed, by assumption, we have an isomorphism
\[
\O_Y \rightarrow F(\O_Y) = \pi_*(E).
\]
By adjunction, we have a morphism
\[
\O_{Y \times Y} \simeq \pi^*\O_Y \rightarrow E.
\]
Let $s$ denote the restriction of this morphism to the diagonal. We
claim that $s$ is the desired isomorphism. First, we prove it is an
isomorphism on fibers. Let
\[
r : E \rightarrow E|_{Y \times y}
\]
denote restriction. By property (LA), we have
\[
\pi_*(E|_{Y \times y}) \simeq \pi_*(E)_y.
\]
By property (LA) and the exactness of affine push-forward,
\[
\pi_*(r) : \pi_*(E) \rightarrow \pi_*(E|_{Y \times y})
\]
is surjective. So, the composition
\begin{equation}\label{eq:suc}
\O_Y \rightarrow \pi_*(E) \overset{\pi_*(r)}{\rightarrow}
\pi_*(E|_{Y \times y})
\end{equation}
is surjective onto the non-zero sheaf $\pi_*(E)_y \simeq \sky_y,$
and in particular is not zero. By naturality of adjunction, the
composition \eqref{eq:suc} is adjoint to the composition
\[
\O_{Y \times Y} \rightarrow E \overset{r}{\rightarrow} E|_{Y \times
y}.
\]
In particular, the latter cannot be zero. By equation \eqref{eq:fib}
and the fact that $E$ is concentrated in one degree, $E|_{Y \times
y}$ has only one non-vanishing fiber,
\[
(E|_{Y \times y})_y \simeq E_{y \times y} \simeq \O_{y \times y}.
\]
We conclude that the composition
\[
\O_{Y \times Y} \rightarrow E \rightarrow E|_{y \times y}
\]
must be surjective for all $y \in Y.$ It follows that $s$ is an
isomorphism on fibers. By Nakayama's lemma, $s$ is surjective. Since
$\Delta$ is reduced, there are no nilpotents among the sections of
$\O_\Delta.$ So, we could detect any non-trivial section in the
kernel of $s$ at some fiber. We conclude that $s$ is an isomorphism.

Finally, we prove that the restriction map
\[
q: E \rightarrow E|_{\Delta}
\]
is an isomorphism. We know that $q$ is surjective. By property (LA)
and exactness of affine push-forward, we know that the map
\begin{equation*}
\O_Y \simeq \pi_*(E) \overset{\pi_*(q)}{\rightarrow}
\pi_*(E|_\Delta)
\end{equation*}
is surjective. By isomorphism \eqref{eq:m1},
\[
\pi_*(E|_\Delta) \simeq \pi_*(\O_\Delta) \simeq \O_Y.
\]
So, $\pi_*(q)$ is a surjective map from $O_Y$ to itself. So,
$\pi_*(q)$ is multiplication by a non-vanishing function, and hence
is an isomorphism. Now, suppose $q$ has a kernel $K.$ Since
$\pi_*(q),$ is an isomorphism, again using property (LA) and
exactness of affine push-forward, we conclude that $\pi_*K$
vanishes. Using property (LA) and the fact that affine push-forward
cannot send a non-trivial sheaf to the zero sheaf, we conclude that
$K$ vanishes. So, $q$ is an isomorphism.

Composing the isomorphism $q$ with the inverse of isomorphism
\eqref{eq:m1}, we have
\[
E \simeq \O_\Delta.
\]
So, $F \simeq \id,$ as claimed.
\end{proof}

We now prove Theorems~\ref{tm:o1} and~\ref{tm:o2}.
\begin{proof}[Proof of Theorem~\ref{tm:o1}]
The auto-equivalence $\D \circ \D'^{op}$ of $\DC(Y)$ satisfies the
hypothesis of Corollary~\ref{cor:maa}. So, we have
\[
\D' \simeq (\D^{-1})^{op} \simeq \D.
\]
\end{proof}
\begin{proof}[Proof of Theorem~\ref{tm:o2}]
The auto-equivalence $\T^{-1} \circ \T'$ of $\DC(Y)$ satisfies the
hypothesis of Corollary~\ref{cor:maa}. The theorem follows.
\end{proof}

\bibliographystyle{plain}
\bibliography{inv12}

\end{document}